\documentclass[a4paper,12pt]{amsart}
\usepackage[]{amscd,amsfonts,amsmath,amsxtra,amssymb,bm}
\usepackage{graphics}
\usepackage{extarrows}
\usepackage[all]{xy}
\usepackage{hyperref}
\usepackage{eucal}
\usepackage{rotating}
\usepackage[T1]{fontenc}
\newtheorem{sub}{}[section]
\newtheorem{subsub}{}[sub]
\usepackage{rotating}
\usepackage[active]{srcltx}
\def\ov#1{\overline{#1}}
\def\codim{\mathop{\rm codim}\nolimits}

\def\Ext{\mathop{\rm Ext}\nolimits}

\def\Pic{\mathop{\rm Pic}\nolimits}

\def\imm{\mathop{\rm im}\nolimits}

\def\rk{\mathop{\rm rk}\nolimits}

\def\spec{\mathop{\rm spec}\nolimits}

\def\Proj{\mathop{\rm Proj}\nolimits}

\def\lra{\longrightarrow}

\def\sigg{\mathop{\hbox{$\displaystyle\sum$}}\limits}

\def\hfl#1#2{\smash{\mathop{\ \hbox to 12mm{\rightarrowfill}}
\limits^{\scriptstyle#1}_{\scriptstyle#2} \ }}
\def\hflb#1#2{\smash{\mathop{\hbox to 12mm{\leftarrowfill}}
\limits^{\scriptstyle#1}_{\scriptstyle#2}}}

\def\m#1{{\hbox{$#1$}}}

\def\ot{\otimes}
\def\og{\leavevmode\raise.3ex\hbox{$\scriptscriptstyle\langle\!\langle$}}
\def\fg{\leavevmode\raise.3ex\hbox{$\scriptscriptstyle\,\rangle\!\rangle$}}

\def\nsp{\lbrace 0\rbrace}

\def\dsp{\displaystyle}
\def\Ssect#1#2{\pagebreak[3]\begin{sub}\label{#2}{\sc\small\small
#1}\rm\medskip}

\def\sepsec{\vskip 1.4cm}
\def\sepsub{\vskip 0.7cm}
\def\sepsubsub{\vskip 0.5cm}
\def\sepprop{\vskip 0.5cm}
\def\sepx{\vskip 0.6cm}
\def\xmat#1{\[\xymatrix{#1}\]}
\def\flinc{\ar@{^{(}->}}
\def\fleq{\ar@{=}}
\def\flon{\ar@{->>}}
\def\fmaps{\ar@{|-{>}}}
\def\fflat{\ar@{-}}
\def\fimpl{\ar@{=>}}
\def\Nligne{\hfil\break}

\def\ED{\vskip 1cm\end{document}}

\def\BS#1{{\boldsymbol{#1}}}

\def\mm{\mathfrak{m}}
\def\deriv#1#2{\frac{\partial #1}{\partial #2}}
\def\wt{\widetilde}

\newcommand{\A}{{\mathbb A}}
\newcommand{\B}{{\mathbb B}}
\newcommand{\Z}{{\mathbb Z}}

\newcommand{\C}{{\mathbb C}}
\newcommand{\Q}{{\mathbb Q}}
\renewcommand{\P}{{\mathbb P}}

\newcommand{\F}{{\mathbb F}}

\newcommand{\Y}{{\mathbb Y}}
\newcommand{\U}{{\mathbb U}}
\newcommand{\Ss}{{\mathbb S}}

\newcommand{\X}{{\mathbb X}}
\renewcommand{\L}{{\mathbb L}}

\def\T{{\mathbb T}}

\newcommand{\kd}{{\mathcal D}}

\newcommand{\ki}{{\mathcal I}}

\newcommand{\ko}{{\mathcal O}}

\newcommand{\kx}{{\mathcal X}}
\newcommand{\ky}{{\mathcal Y}}

\begin{document}

\def\refname{References}
\def\contentsname{Summary}
\def\proofname{Proof}
\def\abstractname{Resume}

\author{Jean--Marc Dr\'{e}zet}
\address{
Institut de Math\'ematiques de Jussieu - Paris Rive Gauche\\
Case 247\\
4 place Jussieu\\
F-75252 Paris, France}
\email{jean-marc.drezet@imj-prg.fr}
\title[{blowing-ups}] {Blowing-ups of primitive multiple schemes}

\begin{abstract}
A primitive multiple scheme is a Cohen-Macaulay scheme $\BS{\kx}$ such that the 
associated reduced scheme $X=\BS{\kx}_{red}$ is smooth, irreducible, and that 
$\BS{\kx}$ can be locally embedded in a smooth variety of dimension 
\m{\dim(X)+1}. If \m{\ki_X} is the ideal sheaf of $X$, \m{L=\ki_X/\ki_X^2} is a 
line bundle on $X$, called the {\em associated line bundle} of $\BS{\kx}$. The 
simplest example is the trivial primitive multiple scheme of multiplicity $n$ 
associated to a line bundle $L$ on $X$: it is the $n$-th infinitesimal 
neighborhood of $X$, embedded in the line bundle $L^*$ by the zero section.

A subscheme $Z$ of $\BS{\kx}$ is called {\em good} if \m{Z_{red}} is smooth and 
connected, and if \m{d=\codim(Z)} and \m{\ki_{Z,\BS{\kx}}} is the ideal sheaf 
of $Z$, then for every closed point \m{z\in Z}, \m{\ki_{Z,\BS{\kx},z}} can be 
generated by $d$ elements.

Two kinds of subschemes $Z$ of $\BS{\kx}$ will be considered: the closed smooth 
subschemes of $X$, seen as subschemes of $\BS{\kx}$, and the good subschemes. 
In the two cases, the blowing-up \m{{\bf B}_{Z,\BS{\kx}}} of $\BS{\kx}$ along 
$Z$ is a primitive multiple scheme of multiplicity $n$, and its underlying 
smooth scheme is the blowing-up \m{{\bf B}_{Z_{red},X}} of $X$ along 
\m{Z_{red}}. Additional results are obtained in the case of hypersurfaces or 
points of $X$.

We treat the case of \m{X=\P_2}, with $Z$ a single point $P$. Let \ 
\m{\pi:\wt{\P}_2\to\P_2} \ be the blowing-up of \m{\P_2} along $P$. We find all 
primitive double schemes $\X$, $\Y$, with \ \m{\Y_{red}=\wt{\P}_2}, 
\m{\X_{red}=\P_2}, such that there is a morphism \m{\Y\to\X} inducing $\pi$ and 
an isomorphism \ \m{\Y\backslash\pi^{-1}(P)\to\X\backslash\{P\}}. We obtain 
in this way the list of all K3-carpets with underlying smooth variety 
\m{\wt{\P}_2}. 
\end{abstract}

\maketitle
\tableofcontents


\vskip 1cm

\section{Introduction}\label{intro}

A {\em primitive multiple scheme} is a Cohen-Macaulay scheme $\BS{\kx}$ over 
$\C$ such that:
\begin{enumerate}
\item[--] $\BS{\kx}_{red}=X$ is a smooth connected variety, 
\item[--] for every closed point \m{x\in X}, there exists a neighborhood $U$ 
of $x$ in $\BS{\kx}$, and a smooth variety $S$ of dimension \ \m{\dim(X)+1} 
such that $U$ is isomorphic to a closed subscheme of $S$.
\end{enumerate}
We call $X$ the {\em support} of $\BS{\kx}$. It may happen that $\BS{\kx}$ is 
quasi-projective, and in this case it is projective if $X$ is.

For every closed subscheme \m{Z\subset \BS{\kx}}, let \m{\ki_Z} (or 
\m{\ki_{Z,\BS{\kx}})} denote the ideal sheaf of $Z$ in $\BS{\kx}$. For every 
positive integer $i$, let \m{\BS{\kx}_i} be the closed subscheme of $\BS{\kx}$ 
corresponding to the ideal sheaf \m{\ki_\BS{\kx}^i}. The smallest integer $n$ 
such that \ \m{\BS{\kx}_n=\BS{\kx}} \ is called the {\em multiplicity} of $Y$. 
For \m{1\leq i\leq n}, \m{\BS{\kx}_i} is a primitive multiple scheme of 
multiplicity $i$, \m{L=\ki_\BS{\kx}/\ki_{\BS{\kx}_2}} is a line bundle on 
$X$, and we have \ \m{\ki_{\BS{\kx}_{i}}/\ki_{\BS{\kx}_{i+1}}=L^i}. 
We call $L$ the line bundle on $X$ {\em associated} to $\BS{\kx}$. If \m{n=2}, 
$\BS{\kx}$ is called a {\em primitive double scheme}.

For any \m{L\in\Pic(X)}, the {\em trivial primitive scheme} of multiplicity 
$n$, with induced smooth variety $X$ and associated line bundle $L$ on $X$ is 
the $n$-th infinitesimal neighborhood of $X$, embedded by the zero section in 
the dual bundle $L^*$, seen as a smooth variety.

The primitive multiple curves where defined in \cite{fe}, \cite{ba_fo}. 
Primitive double schemes were parameterized and studied in \cite{ba_ei} and 
\cite{ei_gr}. More results on primitive multiple curves can be found in 
\cite{dr2}, \cite{dr1}, \cite{dr4}, \cite{ch-ka}, \cite{sa1}, \cite{sa2}, 
\cite{sa3}. Some primitive double schemes are studied in \cite{b_m_r}, 
\cite{ga-go-pu} and \cite{gonz1}. The case of varieties of any dimension, and 
higher multiplicity, is studied in \cite{dr10}, where the following subjects 
were treated:
\begin{enumerate}
\item[--] construction and parametrization of primitive multiple schemes,
\item[--] obstructions to the extension of a vector bundle on \m{\BS{\kx}_m} to 
\m{\BS{\kx}_{m+1}},
\item[--] obstructions to the extension of a primitive multiple scheme of 
multiplicity $n$ to one of multiplicity \m{n+1}.
\end{enumerate}
In \cite{dr11}, the construction and properties of fine moduli spaces of vector 
bundles on primitive multiple schemes are described. In \cite{dr12} other types 
of sheaves are studied.

Let $\BS{\ky}$ be a noetherian scheme. If \ \m{Z\subset \BS{\ky}} \ is a closed 
subscheme, let \ \m{\pi_{Z,\BS{\ky}}:{\bf B}_{Z,\BS{\ky}}\to \BS{\ky}} \ be the 
blowing-up of $\BS{\ky}$ along $Z$.

Let $X$ be a smooth irreducible variety, \m{n\geq 2} be an integer and 
\m{\BS{\kx}} a primitive multiple scheme of multiplicity $n$ such that \ 
\m{\BS{\kx}_{red}=X} \ and with associated line bundle \ \m{L\in\Pic(X)}.

\sepsub

\Ssect{The case of subschemes of $X$}{intro1}

Let \m{Y\subset X} a nonempty closed irreducible smooth subvariety, such that 
\m{Y\not=X}. Let 
\[\wt{Y} \ = \ \pi_{Y,X}^{-1}(Y) \ , \qquad \wt{Y}_n \ = \ 
\pi_{Y,\BS{\kx}}^{-1}(Y) \ . \]
We will prove (cf. \ref{BLU_R})

\sepprop

 \begin{subsub}\label{theoA}{\bf Theorem: } \m{{\bf B}_{Y,\BS{\kx}}} is a 
primitive multiple scheme of multiplicity $n$, \m{({\bf 
B}_{Y,\BS{\kx}})_{red}={\bf B}_{Y,X}}, its associated line bundle is
\[\ko_{{\bf B}_{Y,\BS{\kx}}}(\wt{Y}_n)_{|{\bf B}_{Y,X}}\ot\pi_{Y,X}^*(L) \ = \
\ko_{{\bf B}_{Y,X}}(\wt{Y})\ot\pi_{Y,X}^*(L) \ , \]
The exceptional divisor \m{\wt{Y}_n} is a trivial primitive 
multiple scheme of multiplicity $n$, such that \m{(\wt{Y}_n)_{red}=\wt{Y}}, and 
with associated line bundle \m{\big(\ko_{{\bf 
B}_{Y,X}}(\wt{Y})\ot\pi_{Y,X}^*(L)\big)_{|\wt{Y}}}.
\end{subsub}

The proof follows the construction of blowing-ups given in \cite{ha}, II-2. If 
$Y$ is a hypersurface, we give another proof, inspired by \cite{ba_ei} and 
\cite{dr1}, using the fact that \ \m{{\bf B}_{Y,X}=X}.

\sepprop

\begin{subsub} The case of double schemes -- \rm Suppose that $\BS{\kx}$ is of 
multiplicity 2. Let \m{\BS{\kx}'} another double primitive scheme such that \ 
\m{\BS{\kx}'_{red}=X}, not isomorphic to $\BS{\kx}$. We will prove that if \ 
\m{\codim(Y,X)\geq 2}, then \m{{\bf B}_{Y,\BS{\kx}}} and \m{{\bf 
B}_{Y,\BS{\kx}'}} are not isomorphic (proposition \ref{prop7}).

Suppose that $Y$ is a hypersurface. To $\BS{\kx}$ one can associate \ 
\m{\sigma\in H^1(X,T_X\ot L)}, such that \m{\C.\sigma} characterizes $\BS{\kx}$ 
(cf. \ref{DBC}). We have a canonical map \Nligne \m{\phi_Y:H^1(X,T_X\ot L)\lra 
H^1(X,T_X\ot L\ot\ko_X(Y))} . Using the second proof the preceding theorem for 
hypersurfaces, we prove that the element of \m{H^1(X,T_X\ot L\ot\ko_X(Y))} 
associated to \m{{\bf B}_{Y,\BS{\kx}}} is \m{\phi_Y(\sigma)}.
\end{subsub}

\end{sub}

\sepsub

\Ssect{The case of good subschemes of $\BS{\kx}$}{intro2}

Suppose that \m{\BS{\kx}} is of multiplicity $n$.
Let \m{Z\subset\BS{\kx}} be a good subscheme of dimension $d$. Let \ 
\m{Y=Z_{red}} \ and \ \m{\wt{Z}_n=\pi_{Z,\BS{\kx}}^{-1}(Z)}. We will prove (cf. 
\ref{BLG})

\sepprop

\begin{subsub}\label{theoB}{\bf Theorem: } \m{{\bf B}_{Z,\BS{\kx}}} is a 
primitive multiple scheme of multiplicity $n$, \m{({\bf 
B}_{Z,\BS{\kx}})_{red}={\bf B}_{Y,X}}, and its associated line bundle is \ 
\m{\pi_{Y,X}^*(L)}.

The exceptional divisor \m{\wt{Z}_n} is a primitive multiple scheme of 
multiplicity $n$,\Nligne \m{(\wt{Z}_n)_{red}=\wt{Y}=\pi_{Y,X}^{-1}(Y)}, and 
its associated line bundle is \ \m{\pi_{Y,X}^*(L)_{|\wt{Y}}}.
\end{subsub}

The proof is very similar to that of theorem \ref{theo1}. The two cases are 
related: the blowing-up of \m{{\bf B}_{Z,\BS{\kx}}} along \
\m{\wt{Y}=\pi_{Y,X}^{-1}(Y)\subset{\bf B}_{Y,X}} \ is \m{{\bf B}_{Y,\BS{\kx}}} 
(proposition \ref{prop11}).

\sepprop

{\em Zero-dimensional good subschemes -- } \rm Suppose that $\BS{\kx}$ is of 
multiplicity 2. Let \m{P\in X}. Let \m{\BS{\Gamma}_P} be the set 
of good subschemes \m{Z\subset \BS{\kx}} such that \ \m{Z_{red}=\{P\}}. Then 
\m{\BS{\Gamma}_P} has a natural structure of affine space, with associated 
vector space \m{(T_X\ot L)_P}. Hence, given \ \m{Z,Z'\in\BS{\Gamma}_P}, we can 
define the ``difference'' \ \m{Z-Z'\in(T_X\ot L)_P}.

Let \ \m{\text{ev}_L:\ko_X\ot H^0(X,T_X\ot L)\to T_X\ot L} \ be the evaluation 
morphism, and \Nligne \m{H_{L,P}=\imm(\text{ev}_{L,P})\subset(T_X\ot L)_P}.

In \ref{BL}, we determine when \m{{\bf B}_{Z,\kx}} and \m{{\bf B}_{Z,\kx'}} are 
isomorphic. When $\BS{\kx}$ is not trivial, it is the case if and only if \ 
\m{Z-Z'\in H_{L,P}} (theorem \ref{theo3}).

\end{sub}

\sepsub

\Ssect{Double structures on \m{\P_2}}{dbp}

Let \ \m{P\in\P_2}. Let \ \m{\pi:\wt{\P}_2\to\P_2} \ be the blowing-up of 
\m{\P_2} along $P$, and \ \m{\wt{P}=\pi^{-1}(P)} \ the exceptional divisor.

According to \cite{dr10}, there exists only one non trivial primitive double 
scheme $\X$ such that \ \m{\X_{red}=\P_2}, and its associated line bundle is 
\m{\ko_{\P_2}(-3)}.

Let $\Y$ be a primitive double scheme such that \ \m{\Y_{red}=\P_2},  and with 
associated line bundle \m{\ko_{\P_2}(m)}. We will give the list of all 
primitive double schemes \m{\wt{\Y}} such that
\begin{enumerate}
\item[--] $\wt{\Y}_{red}=\wt{\P}_2$ ,
\item[--] there exists a morphism \ \m{\pi_{\wt{\Y}}:\wt{\Y}\to\Y} \ inducing 
\ \m{\pi_{P,\P_2}:\wt{\P}_2\to \P_2}, and such that its restriction \ 
\m{\wt{\Y}\backslash\wt{P}\to\Y\backslash\{P\}} \ is an isomorphism.
\end{enumerate}
In this case, if $\L$ is the line bundle on \m{\wt{\P}_2} associated to 
$\wt{\Y}$, there is an integer $p$ such that \ 
\m{\L=\pi_{P,\P_2}^*(\ko_{\P_2}(m))\ot\ko_{\wt{\P}_2}(p\wt{P})} . We will prove

\sepprop

\begin{subsub}{\bf Theorem: } We have \ \m{p\geq 0}.

{\bf 1 -- } If \m{p=0}, the schemes $\wt{\Y}$ are the blowing-ups of $\Y$ along 
the good subschemes $Z$ such that \m{Z_{red}=\{P\}}. If \m{m\geq -1}, all 
these schemes are isomorphic. If \m{m\leq -2} they are all distinct.

{\bf 2 -- } If \m{p=1}, these schemes form a family 
\m{(\wt{\Y}_{1,\alpha})_{\alpha\in\C}}, \m{\wt{\Y}_{1,0}} being the blowing-up 
of the schemes of {\bf 1} along $\wt{P}$, and also the blowing-up of $\Y$ along 
$P$. If \m{\Y=\X}, these schemes are all distinct. If $\Y$ is trivial, 
\m{\wt{\Y}_{1,0}} is trivial, and all the \m{\wt{\Y}_{1,\alpha}}, 
\m{\alpha\not=0} are isomorphic and non trivial.

{\bf 3 -- } If \m{p\geq 2}, there is only one such scheme \m{\wt{\Y}_p}. If 
\m{p=2}, \m{\wt{\Y}_p} is the blowing-up of the schemes of {\bf 2} along 
$\wt{P}$. If \m{p>2}, \m{\wt{\Y}_p} is the blowing-up of \m{\wt{\Y}_{p-1}} 
along $\wt{P}$.
\end{subsub}

\sepprop

{\em K3-carpets -- }
We consider the schemes \m{\wt{\X}_{1,\alpha}}, \m{\alpha\in\C}, of theorem 
1.3.1, 2 (when $\Y$ is the non trivial double scheme $\X$). The line 
bundle associated to \m{\wt{\X}_{1,\alpha}} is \ 
\m{\pi_{P,\P_2}^*(\ko_{\P_2}(-3))\ot\ko_{\wt{\P}_2}(\wt{P})=
\omega_{\wt{\P}_2}}. It follows that \m{\wt{\X}_{1,\alpha}} is a {\em 
K3-carpet} (cf. \cite{dr11}, \cite{ga-go-pu}).

\sepprop

\begin{subsub}\label{theoC}{\bf Theorem: } {\bf 1 -- } \m{\wt{\X}_{1,\alpha}} 
is quasi-projective if and only $\alpha\in\Q$ and \m{\alpha>0}.

{\bf 2 -- } If $\alpha\not\in\Q$, 
\m{\Pic(\wt{\X}_{1,\alpha})=\{\ko_{\wt{\X}_{1,\alpha}}\}}. 

{\bf 3 -- }If $\alpha\in\Q$, \m{\Pic(\wt{\X}_{1,\alpha})\simeq\Z}.
\end{subsub}

\sepprop

If instead of $\X$ we take for $\Y$ the trivial double scheme with associated 
line bundle \m{\ko_{\P_2}(-3)}, we obtain two K3-carpets, the trivial scheme 
\m{\wt{\Y}_{1,0}} (the blowing-up of $\Y$ along $P$), and  
\m{\wt{\Y}_{1,\alpha}}, \m{\alpha\not=0}. We have 
\m{\Pic(\wt{\Y}_{1,\alpha})\simeq \Z}, generated by an 
extension of \m{\pi^*(\ko_{\P_2}(1))} to \m{\wt{\Y}_{1,\alpha}}.

\end{sub}

\sepsub

\Ssect{Outline of the paper}{intro-5}

In Chapter 2, we give several preliminary definitions and technical results 
that will be used in the next chapters.

In Chapter 3 we recall the definitions and properties of primitive multiple 
schemes and of good subschemes. We show also how to define morphisms between 
primitive double schemes.

In Chapter 4, the blowing-ups are described, and theorems \ref{theoA} and 
\ref{theoB} are proved.

In Chapter 5, we give supplementary results on blowing-ups of $\BS{\kx}$ along 
zero-dimensional good subschemes, when $\BS{\kx}$ is of multiplicity 2.

In Chapter 6, we give another proof of theorem \ref{theoA}, when $Y$ is an 
hypersurface.

In Chapter 7, we treat the case of \m{X=\P_2}, and prove the results of 
\ref{dbp}.

\end{sub}

\sepprop

{\bf Terminology: } A {\em scheme} is a noetherian separated scheme over $\C$.

\sepsec

\section{Preliminaries}\label{prelim}

\Ssect{Canonical class of a vector bundle}{can_class}

\begin{subsub}\label{CC_def} Definition -- \rm Let $Z$ be a scheme over $\C$ 
and $L$ a line bundle on $Z$. To $L$ one associates an element \m{\nabla_0(L)} 
of \m{H^1(Z,\Omega_Z)}, called the {\em canonical class of $L$}. If $Z$ is 
smooth and projective, and \m{L=\ko_Z(Y)}, where \m{Y\subset Z} is a smooth 
hypersurface, then \m{\nabla_0(L)} is the cohomology class of $Y$.

Let \m{(Z_i)_{i\in I}} be an open cover of $Z$ such that $L$ is defined by a 
cocycle \m{(\theta_{ij})}, \m{\theta_{ij}\in\ko_Z(Z_{ij})^*}. Then 
\m{\dsp\Big(\frac{d\theta_{ij}}{\theta_{ij}}\Big)} is a cocycle which 
represents \m{\nabla_0(L)}.
\end{subsub}

\end{sub}

\sepsub

\Ssect{\v Cech cohomology and trivializations}{cech}

Let $X$ be a scheme over $\C$ and $E$ a vector bundle on $X$, \m{r=\rk(E)}. 
Let \m{(U_i)_{i\in I}} be an open cover of $X$ such that we have 
trivializations:
\xmat{\alpha_i:E_{|U_i}\ar[r]^-\simeq & \ko_{U_i}\ot\C^r .}
Let
\xmat{\alpha_{ij}=\alpha_i\circ\alpha_j^{-1}:\ko_{U_{ij}}\ot\C^r\ar[r]^-\simeq &
\ko_{U_{ij}}\ot\C^r ,}
so that we have the relation \ \m{\alpha_{ij}\alpha_{jk}=\alpha_{ik}}.

\sepprop

\begin{subsub}\label{cech1}\rm Let $n$ be a positive integer, and for every 
sequence \m{(i_0,\ldots,i_n)} of distinct elements of $I$, \m{\sigma_{i_0\cdots 
i_n}\in H^0(U_{i_0\cdots i_n},\ko_X\ot\C^r)}. Let
\[\theta_{i_0\cdots i_n} \ = \ \alpha_{i_0}^{-1}\sigma_{i_0\cdots i_n} \ \in \
H^0(U_{i_0\cdots i_n},E) \ . \]
The family \m{(\theta_{i_0\cdots i_n})} represents an element of \m{H^n(X,E)} 
if the cocycle relations are satisfied: for every sequence  
\m{(i_0,\ldots,i_{n+1})} of distinct elements of $I$, 
\[\sigg_{k=0}^n(-1)^k\theta_{i_0\cdots\widehat{i_k}\cdots i_{n+1}} \ = \ 0 \ , 
\]
which is equivalent to
\[\alpha_{i_0i_1}\sigma_{i_1\cdots i_{n+1}}+\sigg_{k=1}^{n+1}
(-1)^k\sigma_{i_0\cdots\widehat{i_k}\cdots i_{n+1}} \ = \ 0 \ . \]
For \m{n=1}, this gives that elements of \m{H^1(X,E)} are represented by 
families \m{(\sigma_{ij})},\Nligne \m{\sigma_{ij}\in H^0(U_{ij},\ko_X\ot\C^r)}, 
such that
\[\alpha_{ij}\sigma_{jk}+\sigma_{ij}-\sigma_{ik} \ = \ 0 \ . \]
In \v Cech cohomology, it is generally assumed that 
\m{\theta_{ji}=-\theta_{ij}}. This implies that \ 
\m{\sigma_{ji}=-\alpha_{ji}\sigma_{ij}}.
\end{subsub}

\sepprop

\begin{subsub}\label{cech2}\rm Let $F$ be a vector bundle on $X$. Similarly, an 
element of \m{H^n(X,F\ot E)} is represented by a family  \m{(\mu_{i_0\cdots 
i_n})}, with \ \m{\mu_{i_0\cdots i_n}\in H^0(U_{i_0\cdots i_n},F\ot\C^r)}, 
satisfying the relations
\[(I_F\ot\alpha_{i_0i_1})(\mu_{i_1\cdots i_{n+1}})+\sigg_{k=1}^{n+1}
(-1)^k\mu_{i_0\cdots\widehat{i_k}\cdots i_{n+1}} \ = \ 0 \ . \]
The corresponding element of \m{H^0(U_{i_0\cdots i_n},E\ot F)} is \ 
\m{\theta_{i_0\cdots i_n}=(I_F\ot\alpha_{i_0})^{-1}(\mu_{i_0\cdots i_n})}.

For \m{n=1}, this gives that elements of \m{H^1(X,F\ot E)} are represented by 
families \m{(\sigma_{ij})},\Nligne \m{\sigma_{ij}\in H^0(U_{ij},F\ot\C^r)}, 
such that
\[(I_F\ot\alpha_{ij})\sigma_{jk}+\sigma_{ij}-\sigma_{ik} \ = \ 0 \ . \]
If $E$ is a line bundle, then \ \m{\alpha_{ij}\in\C^*}, and the cocycle 
relation is \ \m{\alpha_{ij}\sigma_{jk}+\sigma_{ij}-\sigma_{ik}=0}.

The family \m{(\sigma_{ij})} represents 0 if and only if there exists, for 
every \m{i\in I}, \m{\rho_i\in H^0(U_i,F\ot\C^r)}, such that \ 
\m{\sigma_{ij}=\rho_i-(I_F\ot\alpha_{ij})(\rho_j)}.
\end{subsub}

\sepprop

\begin{subsub}\label{cech3} \rm Let $D$ be a line bundle on $X$, and \ \m{s\in 
H^0(X,D)}. Let \ \m{h\in H^1(X,E\ot F)}, represented (according to \ref{cech2}) 
by a family \m{(h_{ij})}, \m{h_{ij}\in H^0(X,F\ot\C^r)}, satisfying the 
relations \ \m{\alpha_{ij}h_{jk}+h_{ij}-h_{ik}=0}. The section $s$ induces a map
\[\phi_s:H^1(X,E\ot F)\lra H^1(X,E\ot F\ot D) . \]
Let \ \m{\rho_i:D_{|U_i}\to\ko_{U_i}} \ be trivializations, and \ 
\m{s_i=\rho_is\in\C}. Then \m{(\rho_{ij}\alpha_{ij})} is a cocycle that defines 
the vector bundle \m{E\ot D}. We have \ \m{s_ih_{ij}\in H^0(X,F\ot\C^r)}, and 
using the method of \ref{cech2} (with $E$ replaced with \m{E\ot D}), 
\m{\phi_s(h)} is represented by the family \m{(s_ih_{ij})}.
\end{subsub}

\end{sub}

\sepsec

\section{Primitive multiple schemes}\label{PMS}

\Ssect{Definition and construction}{PMS_def}

Let $X$ be a smooth connected variety, and \ \m{d=\dim(X)}. A {\em multiple 
scheme with support $X$} is a Cohen-Macaulay scheme $\BS{\kx}$ such that 
\m{\BS{\kx}_{red}=X}. If $\BS{\kx}$ is quasi-projective we say that it is a 
{\em multiple variety with support $X$}. In this case $\BS{\kx}$ is projective 
if $X$ is.

Let $n$ be the smallest integer such that \m{\BS{\kx}=X^{(n-1)}}, \m{X^{(k-1)}}
being the $k$-th infinitesimal neighborhood of $X$, i.e. \
\m{\ki_{X^{(k-1)}}=\ki_X^{k}} . We have a filtration \ \m{X=X_1\subset
X_2\subset\cdots\subset X_{n}=\BS{\kx}} \ where $X_i$ is the biggest 
Cohen-Macaulay subscheme contained in \m{\BS{\kx}\cap X^{(i-1)}}. We call $n$ 
the {\em multiplicity} of $Y$.

We say that $\BS{\kx}$ is {\em primitive} if, for every closed point $x$ of $X$,
there exists a smooth variety $S$ of dimension \m{d+1}, containing a 
neighborhood of $x$ in $\BS{\kx}$ as a locally closed subscheme. In this case, 
\m{L=\ki_X/\ki_{X_2}} is a line bundle on $X$, \m{X_j} is a primitive multiple 
scheme of multiplicity $j$ and we have \ \m{\ki_{X_j}=\ki_X^j}, 
\m{\ki_{X_{j}}/\ki_{X_{j+1}}=L^j} \ for \m{1\leq j<n}. We call $L$ the line 
bundle on $X$ {\em associated} to $\BS{\kx}$. The ideal sheaf 
\m{\ki_{X,\BS{\kx}}} can be viewed as a line bundle on \m{X_{n-1}}.

Let \m{P\in X}. Then there exist elements \m{y_1,\ldots,y_d}, $t$ of 
\m{m_{S,P}} whose images in \m{m_{S,P}/m_{S,P}^2} form a basis, and such that 
for \m{1\leq i<n} we have \ \m{\ki_{X_i,P}=(t^{i})}. In this case the images of 
\m{y_1,\ldots,y_d} in \m{m_{X,P}/m_{X,P}^2} form a basis of this vector space.

Even if $X$ is projective, we do not assume that $\BS{\kx}$ is projective.

The simplest case is when $\BS{\kx}$ is contained in a smooth variety $S$ of 
dimension \m{d+1}. Suppose that $\BS{\kx}$ has multiplicity $n$. Let \m{P\in X} 
and \m{f\in\ko_{S,P}}  a local equation of $X$. Then we have \ 
\m{\ki_{X_i,P}=(f^{i})} \ for \m{1<j\leq n} in $S$, in particular 
\m{\ki_{\BS{\kx},P}=(f^n)}, and \ \m{L=\ko_X(-X)} .

For any \m{L\in\Pic(X)}, the {\em trivial primitive variety} of multiplicity 
$n$, with induced smooth variety $X$ and associated line bundle $L$ on $X$ is 
the $n$-th infinitesimal neighborhood of $X$, embedded by the zero section in 
the dual bundle $L^*$, seen as a smooth variety.

If \m{U\subset X} is open, we note \m{U^{(n)}} the corresponding open subscheme 
of \m{X_n}.

\sepsubsub

\begin{subsub}\label{PMS-1} Construction of primitive multiple schemes -- \rm
Let $\BS{\kx}$ be a primitive multiple scheme of multiplicity $n$, 
\m{X=\BS{\kx}_{red}}. Let \ \m{{\bf Z}_n=\spec(\C[t]/(t^n))}. Then for every 
closed point \m{P\in X}, there exists an open neighborhood $U$ of $P$ in $X$ 
such that \begin{enumerate}
\item[--] There exists a section \ \m{\ko_X(U)\to\ko_\BS{\kx}(U^{(n)})} \ of 
the restriction map \ \m{\ko_\BS{\kx}(U^{(n)})\to\ko_X(U)},
\item[--] $L_{|U}$ is trivial.
\end{enumerate}
Then there exists a commutative diagram
 \xmat{ & U\flinc[ld]\flinc[rd] \\
U^{(n)}\ar[rr]^-\simeq & & U\times {\mathbf Z}_n}
i.e. $\BS{\kx}$ is locally trivial (\cite{dr1}, th\'eor\`eme 5.2.1, corollaire 
5.2.2).

It follows that we can construct a primitive multiple scheme of multiplicity 
$n$ by taking an open cover \m{(U_i)_{i\in I}} of $X$ and gluing the varieties 
\ \m{U_i\times{\bf Z}_n} (with automorphisms of the \ \m{U_{ij}\times{\bf Z}_n} 
\ leaving \m{U_{ij}} invariant).

Let \m{(U_i)_{i\in I}} be an affine open cover of $X$ such that we have 
trivializations
\xmat{\delta_i:U_i^{(n)}\ar[r]^-\simeq & U_i\times{\bf Z}_n , }
and \ \m{\delta_i^*:\ko_{U_i\times{\bf Z}_n}\to\ko_{U_i^{(n)}}} \ the 
corresponding isomorphism. Let
\xmat{\delta_{ij}=\delta_j\delta_i^{-1}:U_{ij}\times{\bf Z}_n\ar[r]^-\simeq & 
U_{ij}\times{\bf Z}_n \ . }
Then \ \m{\delta_{ij}^*=\delta_i^{*-1}\delta_j^*} \ is an automorphism of \ 
\m{\ko_{U_i\times Z_n}=\ko_X(U_{ij})[t]/(t^n)}, such that for every \ 
\m{\phi\in\ko_X(U_{ij})}, seen as a polynomial in $t$ with coefficients in 
\m{\ko_X(U_{ij})}, the term of degree zero of \m{\delta_{ij}^*(\phi)} is the 
same as the term of degree zero of $\phi$.
\end{subsub}

\sepsubsub

\begin{subsub}\label{I_X} The ideal sheaf of $X$ -- \rm
There exists \ \m{\alpha_{ij}\in\ko_X(U_{ij})[t]/(t^{n-1})} \ such that 
\ \m{\delta_{ij}^*(t)=\alpha_{ij}t}. Let \ 
\m{\alpha^{(0)}_{ij}=\alpha_{ij|X}\in\ko_X(U_i)}.
For every \m{i\in I}, \m{\delta_i^*(t)} is a generator of 
\m{\ki_{X,\BS{\kx}|{U^{(n)}}}}. So we have local trivializations
\[\xymatrix@R=5pt{\lambda_i:\ki_{X,\BS{\kx}|{U_i^{(n-1)}}}\ar[r] & 
\ko_{U_i^{(n-1)}}\\
\delta_i^*(t)\fmaps[r] & 1}\]
Hence \ \m{\lambda_{ij}=\lambda_i\lambda_j^{-1}: 
\ko_{U_{ij}^{(n-1)}}\to\ko_{U_{ij}^{(n-1)}}} \ is the multiplication by 
\m{\delta_j^*(\alpha_{ij})}. It follows that 
\m{(\delta_j^*(\alpha_{ij}))} (resp. 
\m{(\alpha^{(0)}_{ij})})  is a cocycle representing the line bundle 
\m{\ki_{X,\BS{\kx}}} (resp. \m{L}) on \m{X_{n-1}} (resp. $X$). 
\end{subsub}

\end{sub}

\sepsub

\Ssect{The case of double schemes}{DBC}

We suppose that \m{n=2}. Let \ \m{\alpha_i:L_{|U_i}\to\ko_{U_i}} \ be 
isomorphisms such that \ \m{\alpha_{ij}=\alpha_i\circ\alpha_j^{-1}} \ on 
\m{U_{ij}}. Then we have the following description of \m{\delta_{ij}^*}:
\begin{enumerate}
\item[--] there is a derivation \m{D_{ij}} of \m{\ko_X(U_{ij})} such that \ 
\m{\delta_{ij}^*(\beta)=\beta+D_{ij}(\beta)t} \ for every 
\m{\beta\in\ko_X(U_{ij})}.
\item[--] $\delta_{ij}^*(t)=\alpha_{ij}t$ .
\end{enumerate}
The relation \ \m{\delta_{ij}^*\delta_{jk}^*=\delta_{ik}^*} \ is equivalent to 
\m{D_{ij}+\alpha_{ij}D_{jk}=D_{ik}}, i.e. \m{(D_{ij})} is a cocyle in the sense 
of \ref{cech2}. We can view \m{D_{ij}} as section of 
\m{T_{X|U_{ij}}}, the family \m{(\alpha_i^{-1}\ot D_{ij})} represents an 
element $\lambda$ of \m{H^1(X,T_X\ot L)}, and \m{\C\lambda} is independent of 
the choice of the \m{\delta_{ij}^*} and \m{\alpha_i}. We will note \ 
\m{\C\lambda=\sigma(X_2)\in\P(H^1(X,T_X\ot L))\cup\nsp}.

Let $\BS{\kx}$, \m{\BS{\kx}'} primitive double schemes, 
\m{\phi:\BS{\kx}_{red}\to X}, \m{\phi':\BS{\kx}'_{red}\to X} \ isomorphisms. We 
say that $\BS{\kx}$, \m{\BS{\kx}'} are are {\em isomorphic} if there is an 
isomorphism \ \m{\Phi:\BS{\kx}\to\BS{\kx}'} \ such that \ 
\m{\phi=\Phi_{red}\circ\phi'}. In this case the associated line bundles on $X$ 
of $\BS{\kx}$, \m{\BS{\kx}'} are isomorphic. Let \m{{\bf DS}(X,L)} be the set 
of isomorphism classes of primitive double schemes $\BS{\kx}$ such that \ 
\m{\BS{\kx}_{red}\simeq X} \ and with associated line bundle $L$.

According to \cite{dr10} and \cite{ba_ei}, two primitive double 
schemes \m{X_2}, \m{X'_2}, with underlying smooth variety $X$ and associated 
line bundle $L$ are isomorphic (over $X$) if and only if 
\m{\sigma(X_2)=\sigma(X'_2)}. So there is a canonical bijection \ \m{{\bf 
DS}(X,L)\simeq\P(H^1(X,T_X\ot L))\cup\nsp}.

Another family \m{(D'_{ij})} of derivations defines (with \m{(\alpha_{ij})}) 
the same double scheme as \m{(D_{ij})} if and only if there exist 
\m{\tau\in\C^*}, and for every $i$ a derivation \m{T_i} of \m{\ko_X(U_i)} such 
that \ \m{D'_{ij}=\tau D_{ij}+T_i-\alpha_{ij}T_j}.

\sepprop

\begin{subsub}\label{ef_ch} Change of trivialization -- \rm Let \m{i_0\in I} 
and \m{\theta} and automorphism of \m{U_{i_0}\times{\bf Z}_2}, defined by a 
derivation $D$ of \m{\ko_X(U_{i_0})}, i.e \ \m{\theta^*(u)=u+D(u)t} \ for every 
\m{u\in\ko_X(U_{i_0})}, and \m{\theta^*(t)=t}. Now let \ \m{\mu_i=\delta_i} if 
\m{i\in I,i\not=i_0}, and \ \m{\mu_{i_0}=\theta\circ\delta_{i_0}}. We have \ 
\m{\mu_{jk}^*=\delta_{jk}^*} \ if \m{j,k\in I}, \m{j,k\not=i_0}, and \ 
\m{\mu_{i_0k}^*=\theta^{*-1}\delta_{i_0k}^*}, 
\m{\mu_{ji_0}^*=\delta_{ji_0}^*\theta^*}. With these new trivializations, 
\m{X_2} is defined by the same cocycle \m{(\alpha_{ij})}, and the family 
\m{E_{ij}} of derivations, with \ \m{E_{jk}=D_{jk}} \ if \m{j,k\in I}, 
\m{j,k\not=i_0},
\[E_{i_0k} \ = \ D_{i_0k}-D \ , \qquad E_{ji_0} \ = \ D_{ji_0}+\alpha_{ji_0}D \ 
. \]
\end{subsub}

\end{sub}

\sepsub

\Ssect{Extension of line bundles in higher multiplicity}{ext_hig}

Suppose that \m{X_n} can be extended to a primitive multiple scheme \m{X_{n+1}} 
of multiplicity \m{n+1}. We have an exact sequence of coherent sheaves on 
\m{X_n}
\[0\lra L^n\lra\Omega_{X_{n+1|X_n}}\lra\Omega_{X_n}\lra 0 \ , \]
corresponding to \ \m{\sigma_n\in\Ext^1_{\ko_{X_n}}(\Omega_{X_n},L^n)}. Let 
$\F$ be a line bundle on \m{X_n}. Let \Nligne \m{\nabla_0(\F)\in 
H^1_{\ko_{X_n}}(X_n,\Omega_{X_n})} \ be the canonical class of $\F$ (cf. 
\ref{CC_def}). Then we have 
\[\sigma_n\nabla_0(\F) \ \in \ \Ext^2_{\ko_{X_n}}(\F,\F\ot L^n) \ =
 H^2(X,L^n) \ . \]
From \cite{dr10}, theorem 7.1.2, $\F$ can be extended to a line bundle on 
\m{X_{n+1}} if and only if \ \m{\sigma_n\nabla_0(\F)=0}.

\end{sub}

\sepsub

\Ssect{Good subschemes}{GS}

Let $Z$ a subscheme of \m{\BS{\kx}=X_n}. We say that $Z$ is {\em good} if 
\m{Z_{red}} is smooth and connected, and if \m{d=\codim(Z)} and 
\m{\ki_{Z,\BS{\kx}}} is the ideal sheaf of $Z$, then for every closed point 
\m{z\in Z}, \m{\ki_{Z,\BS{\kx},z}} can be generated by $d$ elements.

\sepprop

\begin{subsub}\label{prop6}{\bf Proposition: } Let $Z$ be a subscheme of 
$\BS{\kx}$ such that if \m{Z_{red}} is smooth and connected of dimension $d$. 
Then $Z$ is good if and only if, for every closed point \m{z\in Z}, there 
exists an open neighborhood \m{U\subset X} of $z$ such that there is a 
trivialization \ \m{\rho:U^{(n)}\to U\times{\bf Z}_n}, and generators 
\m{y_1,\ldots,y_d} of \ \m{\ki_{Z_{red},X}(U)}, which are also generators of 
the ideal sheaf of \m{\rho(Z\cap U^{(n)})}.
\end{subsub}
\begin{proof} If the condition is satisfied, then $Z$ is good. Conversely, 
suppose that $Z$ is good. There exists an affine neighborhood \m{U\subset X} of 
$z$ such that 
\begin{enumerate}
\item[--] $\ki_{Z,\BS{\kx}}(U)$ is generated by $d$ elements 
\m{z_1,\ldots,z_d}. 
\item[--] The restriction of the tangent bundle of $X$ on $U$ is free. 
\item[--] There exists a trivialization \ \m{\theta:U^{(n)}\to U\times{\bf 
Z}_n}.
\end{enumerate}
We can write, for \m{1\leq i\leq d}, \m{\theta^{*-1}(z_i)=y_i+a_it}, with \ 
\m{y_i\in\ki_{Z_{red},X}(U)} \ \m{a_i\in\ko_X(U)[t]/(t^{n-1})}. Then 
\m{y_1,\ldots,y_d} generate \m{\ki_{Z_{red},X}(U)}, and for every \m{z'\in 
Z_{red}\cap U}, the images of \m{y_1,\ldots,y_d} in \m{m_{X,z'}/m_{X,z'}^2} are 
linearly independent (where \m{m_{X,z'}\subset\ko_{X,z'}} is the maximal 
ideal). According to \cite{dr10}, 5, there exists an automorphism $\phi$ of \ 
\m{U\times{\bf Z}_n} leaving $U$ invariant and such that \ 
\m{\phi^*(y_i)=y_i+a_it} \ for \m{1\leq i\leq d}. We can then take \ 
\m{\rho=\theta\circ\phi}.
\end{proof}

\sepprop

{\bf Remarks:} {\bf 1 -- } It is necessary to consider only smooth subvarieties 
\m{Z_{red}\subset X} in proposition \ref{prop6}: suppose that \ 
\m{\BS{\kx}=\C\times{\bf Z}_2=\m{\spec(\C[X_0])\times{\bf Z}_2}}, and let \ 
\m{f\in\subset\C[X_0]} \ be a non zero polynomial multiple of \m{X_0^2}. Let \ 
\m{Z\subset\BS{\kx}} \ be the subscheme defined by the ideal \ 
\m{(f+t)\subset\C[X_0][t]/(t^2)}. Let \ \m{U\subset\C} \ be a neighborhood of 0 
and $\rho$ an automorphism of \m{U\times{\bf Z}_2} leaving $U$ invariant. Then, 
according to \ref{DBC}, there exists a derivation $D$ of \m{\ko_\C(U)} such 
that, for every \m{\alpha\in\ko_\C(U)} we have \ 
\m{\rho^*(\alpha)=\alpha+D(\alpha)t}. Hence \ \m{\rho^*(f)=f+D(f)t}. Since 
\m{D(f)} vanishes at 0, it is impossible to obtain \ \m{\rho^*(f)=f+t}. Hence 
proposition \ref{prop6} is false for \m{Z_{red}\subset\C}.

{\bf 2 -- } If \ \m{d=\dim(X)-1}, $Z$ is a Cartier divisor, and 
\m{\ki_{Z,\BS{\kx}}} is a locally free sheaf. If \m{\BS{\kx}} is the non 
projective double scheme with smooth underlying variety \m{\P_2}, such 
varieties $Z$ of dimension 1 do not exist, since \m{\ko_{\BS{\kx}}} is the only 
line bundle on \m{\BS{\kx}}.

\end{sub}

\sepsub

\Ssect{Morphisms of primitive double schemes}{mor_db}

\begin{subsub}\label{f_der} Derivations -- \rm Let $A$, $B$ be commutative 
$\C$-algebras, and \ \m{f:A\to B} \ a morphism. A map \ \m{\theta:A\to B} \ is 
called a {\em $f$-derivation} if
\begin{enumerate}
\item[--] For every \m{\epsilon\in\C}, \m{\theta(\epsilon)=0}.
\item[--] For every \m{a,a'\in A}, \m{\theta(a+a')=\theta(a)+\theta(a')} .
\item[--] For every \m{a,a'\in A}, 
\m{\theta(a.a')=f(a)\theta(a')+f(a')\theta(a)} .
\end{enumerate}
\end{subsub}

\sepprop

\begin{subsub}\label{mor_db2} Morphisms of primitive double rings -- \rm Let 
$A$, $B$ be commutative $\C$-algebras, and  \m{A[2]=A[t]/(t^2)},  
\m{B[2]=B[t]/(t^2)}. Let \ \m{\phi:A\to B}, \m{\Phi:A[2]\to B[2]} \ be 
morphisms of $\C$-algebras such that \ \m{\Phi((t))\subset(t)} \ and $\phi$ is 
induced by $\Phi$. We have
\[\Phi(\alpha) \ = \ \phi(\alpha)+\gamma(\alpha)t \quad \text{for every \ } 
\alpha\in A \ , \]
\[\Phi(t) \ = \ \tau.t , \]
where \ \m{\gamma:A\to B} \ is a $\phi$-derivation and \ \m{\tau\in B}.
\end{subsub}

\sepprop

\begin{subsub}\label{mor_db3} Morphisms of primitive double schemes -- \rm 
Let $X$, $Y$ be smooth varieties, \m{X_0,X_1} (resp. \m{Y_0,Y_1}) open 
subsets of $X$ (resp. $Y$) such that \ \m{X=X_0\cup X_1} (resp. \m{Y=Y_0\cup 
Y_1}). Let $\BS{\kx}$, $\BS{\ky}$ be primitive double schemes such that \ 
\m{\BS{\kx}_{red}=X} \ and \ \m{\BS{\ky}_{red}=Y}. Assume that we have 
trivializations
\[\lambda_i:X_i^{(2)}\lra X_i\times{\bf Z}_2 \ , \qquad
\mu_i:Y_i^{(2)}\lra Y_i\times{\bf Z}_2 \ , \]
for \m{i=1,2}, i.e. isomorphisms
\[\lambda_i^*:\ko_X(X_i)[t]/(t^2)\lra\ko_\X(X_i) \ , \qquad
\mu_i^*:\ko_Y(Y_i)[t]/(t^2)\lra\ko_\Y(Y_i) \ . \]
Let \ 
\[\lambda=\lambda_1\lambda_0^{-1}:X_{01}\times{\bf Z}_2\to X_{01}\times{\bf 
Z}_2 \ , \quad \lambda^*=\lambda_0^{*-1}\lambda_1^*:\ko_X(X_{01})[t]/(t^2)\to
\ko_X(X_{01})[t]/(t^2) \ , \]
\[\mu=\mu_1\mu_0^{-1}:Y_{01}\times{\bf Z}_2\to Y_{01}\times{\bf Z}_2 \ 
, \quad 
\mu^*=\mu_0^{*-1}\mu_1^*:\ko_Y(Y_{01})[t]/(t^2)\to\ko_Y(Y_{01})[t]/(t^2) \ , 
\]
We have
\[\lambda^*(a) \ = \ a+D_X(a)t \quad \text{for every } a\in\ko_X(X_{01}) \ , 
\quad \lambda^*(t)=\alpha t \ , \]
\[\mu^*(b) \ = \ b+D_Y(b)t \quad \text{for every } b\in\ko_Y(Y_{01}) \ , 
\quad \mu^*(t)=\beta t \ , \]
where \m{D_X} (resp. \m{D_Y}) is a derivation of \m{\ko_X(X_{01})} (resp.
\m{\ko_Y(Y_{01})}), \m{\alpha\in\ko_X(X_{01}^*)} \ and \ 
\m{\beta\in\ko_Y(Y_{01}^*)}.

Let \ \m{\phi:Y\to X} \ be a morphism such that \ \m{\phi(Y_0)\subset X_0} \ 
and \ \m{\phi(Y_1)\subset X_1}. Let \ \m{\Delta_i:Y_i^{(2)}\to X_i^{(2)}} \ be 
a morphism inducing $\phi$. Let
\[\Theta_i=\lambda_i\Delta_i\mu_i^{-1}:Y_i\times{\bf Z}_2\lra X_i\times{\bf 
Z}_2 
\ . \]
Then we have
\[\Theta_i^*(a) \ = \ \phi^*(a)+\gamma_i(a)t \quad \text{for every  } 
a\in\ko_X(X_i) \ , \quad \Theta_i^*(t)=\tau_i.t \ , \]
where \ \m{\gamma_i:\ko_X(X_i)\to\ko_Y(Y_i)} \ is a \m{\phi^*}-derivation 
and \ \m{\tau_i\in\ko_Y(Y_i)^*} .

We have \ \m{\Delta_0=\Delta_1} \ on \m{Y_{01}^{(2)}} if and only if \ 
\m{\Theta_0^*\lambda^*=\mu^*\Theta_1^*} \ on \ \m{X_{01}\times{\bf Z}_2}. This 
is equivalent to
\[\gamma_0+\tau_0\phi^*D_X \ = \ D_Y\phi^*+\beta\gamma_1 \quad \text{on } 
X_{01} \ \text{and } \ \phi^*(\alpha)\tau_0 \ = \ \beta\tau_1 \ . \]

Hence it is equivalent to find an extension \m{\BS{\ky}\to\BS{\kx}} of $\phi$, 
and to find \m{\phi^*}-derivations \m{\gamma_0}, \m{\gamma_1}, 
\m{\tau_i\in\ko_Y(Y_i)^*} \ satisfying the preceding equations. The first 
equation concerns \m{\phi^*}-derivations on \m{X_{01}}, one has to verify that 
\m{\gamma_i} sends \m{\ko_X(X_i)} to \m{\ko_Y(Y_i)}.

Suppose now that \m{X_0=Y_0}, \m{\phi_{|Y_0}=I_{Y_0}} \ and \ 
\m{\Delta_0=I_{Y_0^{(2)}}}. We have then \ \m{\gamma_0=0} and \m{\tau_0=1}. The 
preceding equations become
\begin{equation}\label{equ30}\beta\gamma_1 \ = \ D_X-D_Y \quad \text{on } 
X_{01} \ \text{and } \ \alpha \ = \ \beta\tau_1 \ . \end{equation}

\end{subsub}

\end{sub}

\sepsec

\section{Blowing-up of subschemes of primitive multiple schemes}\label{BUSCH}

(cf. \cite{ha}, II-2)

Let $X$ be a smooth irreducible variety, \m{n\geq 2} be an integer and \m{X_n} 
a primitive multiple scheme of multiplicity $n$ such that \ \m{(X_n)_{red}=X} \ 
and with associated line bundle \ \m{L\in\Pic(X)}.

If \ \m{Z\subset X_n} \ is a closed subscheme, let \ \m{\pi_{Z,X_n}:{\bf 
B}_{Z,X_n}\to X_n} \ be the blowing-up of \m{X_n} along $Z$, and \ 
\m{\wt{Z}_n=\pi_{Z,X_n}^{-1}(Z)}. 

\sepsub

\Ssect{Blowing-ups of reduced subschemes}{BLU_R}

Let \m{Y\subset X} be a smooth closed subvariety. Let \ \m{\pi_{Y,X}:{\bf 
B}_{Y,X}\to X} \ be the blowing-up of $X$ along $Y$ and \ 
\m{\wt{Y}=\pi_{Y,X}^{-1}(Y)}. We can view $Y$ as a closed subscheme of \m{X_n}. 
So we have also \ \m{\pi_{Y,X_n}:{\bf B}_{Y,X_n}\to X_n}, the blowing-up of 
\m{X_n} along $Y$, and \ \m{\wt{Y}_n=\pi_{Y,X_n}^{-1}(Y)}. We will prove

\sepprop

\begin{subsub}\label{theo1}{\bf Theorem: } \m{{\bf B}_{Y,X_n}} is a primitive 
multiple scheme of multiplicity $n$, \Nligne\m{({\bf B}_{Y,X_n})_{red}={\bf 
B}_{Y,X}}, its associated line bundle is
\[\ko_{{\bf B}_{Y,X_n}}(\wt{Y}_n)_{|{\bf B}_{Y,X}}\ot\pi_{Y,X}^*(L) \ = \
\ko_{{\bf B}_{Y,X}}(\wt{Y})\ot\pi_{Y,X}^*(L) \ , \]
The exceptional divisor \m{\wt{Y}_n} is a primitive 
multiple scheme of multiplicity $n$, such that \m{(\wt{Y}_n)_{red}=\wt{Y}}, and 
with associated line bundle \m{\big(\ko_{{\bf 
B}_{Y,X}}(\wt{Y})\ot\pi_{Y,X}^*(L)\big)_{|\wt{Y}}}.
\end{subsub}
\begin{proof}
Let \ \m{U=\spec(A)} \ be an affine irreducible smooth variety, \m{V\subset 
U} a closed subvariety and \ \m{\ki_V=\ki_{V,U}} its ideal sheaf. Let \ 
\m{U_n=U\times{\bf Z}_n}. We have \ \m{V\subset U\subset U_n} \ and \ 
\m{\ko_{U_n}(U)=A[t]/(t^n)}. We will first prove the theorem when \ \m{X=U}, 
\m{Y=V} \ and \ \m{X_n=U_n}, i.e.: \m{{\bf B}_{V,U_n}} is the trivial primitive 
multiple scheme of multiplicity $n$ such that \m{({\bf B}_{V,U_n})_{red}={\bf 
B}_{V,U}}, and with associated line bundle \m{\ko_{{\bf 
B}_{V,U_n}}(\wt{V})_{|{\bf B}_{V,X}}}.

\sepx

{\bf Step 1 }{\em --  Description of \m{{\bf B}_{V,U}} -- } We have \ \m{{\bf 
B}_{V,U}=\Proj(S)}, where \ \m{S=S(V,U)} \ is the graduate ring
\[S \ = \ \bigoplus_{p\geq 0}\ki_V(U)^p \ . \]

The ideal sheaf \m{\ki_{U,U_n}} of $U$ in \m{U_n} is generated by $t$, and we 
have
\[\ki_{V,U_n} \ = \ \ki_V+\ki_{U,U_n} \ . \]
We have \ \m{{\bf B}_{V,U_n}=\Proj(\Ss_n)}, where \ \m{\Ss_n=\Ss_n(V,U)} \ is 
the graduate ring
\[\Ss_n \ = \ \bigoplus_{p\geq 0}\ki_{V,U_n}(U)^p \ , \]
and
\[\ki_{V,U_n}(U)^p \ = \ \bigoplus_{\max(p-n,0)\leq k\leq p}\ki_V(U)^kt^{p-k} \ 
. \]

Let \ \m{f\in\ki_V(U)}, \m{f\not=0}. Let \m{S_f} be the corresponding localized 
ring, and \ \m{S_{(f)}\subset S_f} \ the subring of elements of degree 0. Then 
\ \m{D^S_+(f)=\{\mathfrak{p}\in\Proj(S);f\not\in\mathfrak{p}\}} \ is an affine 
open subset of \m{{\bf B}_{V,U}}, and \ \m{D^S_+(f)\simeq\spec(S_{(f)})}. We 
have 
\[S_{(f)} \ = \ \sum_{p\geq 0}S_{(f),p} \ , \quad \text{with} \ \ \
S_{(f),p} \ = \ \Big\{\frac{\sigma_p}{f^p}\ ;\ \sigma_p\in\ki_V(U)^p\Big\} \ . 
\]
The ideal \m{I_{V,p}\subset S_{(f)}} of \m{V\cap D^S_+(f)} is
\[I_{V,f} \ = \ \sum_{p\geq 0}S'_{(f),p} \ , \quad \text{with} \ \ \
S'_{(f),p} \ = \ \Big\{\frac{\mu_p}{f^p}\ ;\ \mu_p\in\ki_V(U)^{p+1}\Big\} \ . \]
We have an isomorphism \ \m{S_{(f)}\to I_{V,f}} \ , which sends $\alpha$ to 
\m{f\alpha}. Let \m{f_1,\ldots,f_p} be generators of \m{\ki_V(U)}. Then 
\m{(D_+^S(f_i))_{1\leq i\leq p}} is an open cover of \m{{\bf B}_{V,U}}, and 
with respect to this cover, \m{\Big(\dsp\frac{f_j}{f_i}\Big)} is a cocycle that 
defines the line bundle \m{\ko_{{\bf B}_{V,U}}(-\wt{V})}.

\sepx

{\bf Step 2 }{\em --  Description of \m{{\bf B}_{V,U_n}} -- } Similarly, we 
have \ \m{f\in\ki_{V,U_n}(U)}, and we can define the affine open 
subset \m{D^{\Ss_n}_+(f)} of \m{{\bf B_{V,U_n}}}. Let 
\m{\Ss_{n,f}} be the corresponding localized ring, and \ 
\m{\Ss_{n,(f)}\subset\Ss_{n,f}} \ the subring of elements of degree 0. Then \ 
\m{D^{\Ss_n}_+(f)\simeq\spec(\Ss_{n,(f)})}. We have 
\[\Ss_{n,(f)} \ = \ \sum_{p\geq 0}\Ss_{n,(f),p} \ , \quad \text{with} \ \ \
\Ss_{n,(f),p} \ = \ \Big\{\frac{\theta_p}{f^p}\ ;\ 
\theta_p\in\ki_{V,U_n}(U)^p\Big\} \ . \]
It follows that
\[\Ss_{n,(f)} \ = \ S_{(f)}[z]/(z^n) \ , \quad \text{with} \ \ \ z \ = \ 
\frac{t}{f} \ . \]
We can do the same if \ \m{f\in\ki_V(W)}, where \m{W\subset U} is an open 
subset. We obtain an isomorphism \ \m{D^{\Ss_n}_+(f)\simeq D^S_+(f)\times{\bf 
Z}_n} \ inducing the identity on \m{D^S_+(f)}.

\sepx

{\bf Step 3 }{\em -- The primitive multiple scheme structure on \m{{\bf 
B}_{V,U_n}} -- } Let \m{f_1,\ldots,f_p} be generators of \m{\ki_V(U)}. Then 
\m{(D_+^S(f_i))_{1\leq i\leq p}} (resp. \m{(D_+^{\Ss_n}(f_i))_{1\leq i\leq p}}) 
is an open cover of \m{{\bf B}_{V,U}} (resp. \m{{\bf B}_{V,U_n}}). We define 
the primitive scheme structure of \m{{\bf B}_{V,U_n}} by taking the isomorphisms
\[\delta(f_i):D^{\Ss_n}_+(f_i)\lra D^S_+(f_i)\times{\bf Z_n}\]
such that \ \m{\delta(f_i)^*:S_{(f_i)}[z]/(z^n)\lra \Ss_{n,(f_i)}} \ 
\m{\delta(f_i)^*(\alpha)=\alpha} \ for every \ \m{\alpha\in S_{(f_i)}} \ and \ 
\m{\delta(f_i)^*(z)=\dsp\frac{t}{f_i}}. It follows that, if \ 
\m{\delta_{ij}=\delta(f_j)\circ\delta(f_i)^{-1}}, we have \
\m{\delta_{ij}^*(\alpha)=\alpha} for every \ \m{\alpha\in S_{(f_i)}\cap 
S_{(f_j)}}, and \ \m{\delta_{ij}^*(z)=\dsp\frac{f_i}{f_j}z}. It is easily 
verified that \ \m{\delta_{ij}^*\delta_{jk}^*=\delta_{ik}^*}. This proves that 
\m{{\bf B}_{V,U_n}} is a trivial primitive multiple scheme of multiplicity $n$ 
such that \m{({\bf B}_{V,U_n})_{red}={\bf B}_{V,U}}. Moreover we have a 
canonical obvious commutative square
\xmat{{\bf B}_{V,U}\flinc[d]\ar[rr]^-{\pi_{V,U}} & & U\flinc[d]\\
{\bf B}_{V,U_n}\ar[rr]^{\pi_{V,U_n}} & & U_n}

\sepx

{\bf Step 4 }{\em -- The ideal sheaf \m{\ki_{{\bf B}_{V,U}}} of \m{{\bf 
B}_{V,U}} in \m{{\bf B}_{V,U_n}} -- } On \m{D_+^{\Ss_n}(f)}, according to 
step 1, we have \ \m{\ki_{{\bf B}_{V,X}}=\sum_{j\geq 1}I_j}, with 
\[I_j \ = \ \Big\{\frac{t\sigma_j}{f^j}\ ;\ \sigma_j\in\ki_V(U)^{j-1}\Big\} \ . 
\]
It follows that we have a well defined isomorphism \ 
\m{\ko_{D_+^{\Ss_{n-1}}(f)}\to\ki_{{\bf B}_{V,X}|D_+^{\Ss_n}(f)}}, sending 
\m{\dsp\frac{\sigma_j}{f^j}} to \m{\dsp\frac{t\sigma_j}{f^{j+1}}}.
In particular, we have trivializations \ 
\m{\alpha_i:\ko_{D_+^{\Ss_{n-1}}(f_i)}\to\ki_{{\bf B}_{V,X}|D_+^{\Ss_n}(f_i)}} 
\ for \m{1\leq i\leq p}. It follows that the automorphism \m{\alpha_{ij}^*} of
\m{\ko(D_+^{\Ss_{n-1}}(f_i)\cap D_+^{\Ss_{n-1}}(f_j))} is the multiplication 
by \m{\dsp\frac{f_i}{f_j}}. According to step 1, it follows that the line 
bundle on \m{{\bf B}_{V,U}} associated to 
\m{{\bf B}_{V,U_n}} is \m{\ko_{{\bf B}_{V,U}}(\wt{V})}.

\sepx

{\bf Step 5 }{\em -- Automorphisms -- } Let $\BS{\Phi}$ be an automorphism of \ 
\m{U\times{\bf Z}_n} \ leaving $U$ invariant. Then there exists \ 
\m{\alpha\in A[t]/(t^{n-1})}, invertible, such that \ \m{\BS{\Phi}^*(t)=\alpha 
t}. Let $\wt{\BS{\Phi}}$ (resp. \m{\wt{\BS{\Phi}}_n}) be the automorphism of 
\m{{\bf B}_{V,U}} (resp. \m{{\bf B}_{V,U_n}}) induced by $\BS{\Phi}$. Then \ 
\m{\BS{\Phi}^*(\ki_V(U))\subset\ki_V(U)}, and \ \m{\BS{\Phi}^*(f)\in\Ss_{n,f}} 
\ is an invertible element of degree 1.

We have \ \m{\wt{\BS{\Phi}}(D_+^S(f))=D_+^S(f)} \ and \ 
\m{\wt{\BS{\Phi}}_n(D_+^{\Ss_n}(f))=D_+^{\Ss_n}(f)}. Let \m{\wt{\BS{\Phi}}^f} 
(resp. \m{\wt{\BS{\Phi}}_n^f}) be the corresponding automorphism of 
\m{D_+^S(f)} (resp. \m{D_+^{\Ss_n}(f)}).

Let \m{\dsp\frac{\sigma_p}{f^p}\in\Ss_{n,(f)}}. We have
\[\wt{\BS{\Phi}}_n^{f*}\Big(\frac{\sigma_p}{f^p}\Big) \ = \ 
\frac{\BS{\Phi}^*(\sigma_p)}{\BS{\Phi}^*(f)^p}\]
and \ \m{\wt{\BS{\Phi}}_n^{f*}(t)=\alpha t}.

\sepx

{\bf Step 6 }{\em -- The primitive multiple scheme structure on \m{{\bf 
B}_{Y,X_n}} -- } Let \m{(U_i)_{i\in I}} be a finite open affine cover of $X$ 
such that for every \m{i\in I} there is an isomorphism \ 
\m{\lambda_i:U_i^{(n)}\to U_i\times{\bf Z}_n} \ leaving \m{U_i} invariant. We 
have then \ \m{\lambda_{ik}^*(t)=\alpha_{ik}t}, where \  
\m{\alpha_{ik}\in\ko_X(U_{ik})[t]/(t^{n-1})} \ is invertible. The line bundle \ 
\m{L\in\Pic(X)} \ associated to \m{X_n} is defined by the cocycle 
\m{(\alpha_{ik|U_{ik}})}.

Let \ \m{\wt{\lambda_i}:{\bf B}_{Y\cap U_i,U_i^{(n)}}\to{\bf B}_{Y\cap 
U_i,U_i\times{\bf Z}_n}} \ be the induced isomorphism. We have \Nligne 
\m{\pi_{Y\cap U_i,U_i^{(n)}}=\pi_{Y\cap U_i,U_i\times{\bf 
Z}_n}\circ\wt{\lambda_i}}. Then from step 2
\begin{enumerate}
\item[--] ${\bf B}_{Y\cap U_i,U_i^{(n)}}$ is a primitive multiple 
scheme of multiplicity $n$, with associated smooth scheme \m{{\bf B}_{Y\cap 
U_i,U_i}}.
\item[--] We have \ \m{{\bf B}_{Y\cap U_i,U_i}=\pi_{Y,X}^{-1}(U_i)} \ and 
\m{\pi_{Y\cap U_i,U_i}} is the restriction of \m{\pi_{X,Y}}. We have a similar 
statement for \m{{\bf B}_{Y\cap U_i,U_i^{(n)}}}.
\end{enumerate}
Hence \m{{\bf B}_{Y,X_n}} is a primitive multiple scheme of multiplicity $n$, 
and \ \m{({\bf B}_{Y,X_n})_{red}={\bf B}_{Y,X}}.

\sepprop

{\em The open cover -- } More precisely, for every \m{i\in I}, let \ \m{{\bf 
S}_i=S(Y\cap 
U_i,U_i)}, \m{\BS{\Sigma}_i=\Ss_n(Y\cap U_i,U_i)} \ and \m{({\bf 
f}_{ij})_{j\in J_i}} a finite sequence of generators of \m{\ki_{Y\cap 
U_i,U_i}(U_i)}. Let
\[T_{ij} \ = \ D_+^{{\bf S}_i}({\bf f}_{ij}) \ \subset \ {\bf 
B}_{Y\cap U_i,U_i} \ , \quad \T_{ij} \ = \ D_+^{\BS{\Sigma}_i}({\bf 
f}_{ij}) \ \subset \ {\bf B}_{Y\cap U_i,U_i\times{\bf Z}_n} \ , \]  
and \ \m{Z_{ij}=\wt{\lambda}_i^{-1}(\T_{ij})\subset {\bf B}_{Y,X_n}}. We have \ 
\m{(\T_{ij})_{red}=(Z_{ij})_{red}=T_{ij}} \ , and \m{(T_{ij})_{i\in I,j\in 
J_i}} is an open cover of \m{{\bf B}_{Y,X}}.

\sepprop

{\em The trivializations -- } We have seen in step 3 that
\[\delta({\bf f}_{ij}) : \T_{ij}\lra T_{ij}\times{\bf Z}_n\]
is an isomorphism. Then 
\[\rho_{ij} \ = \ \delta({\bf f}_{ij})\circ\wt{\lambda_i}:Z_{ij}\lra 
T_{ij}\times{\bf Z}_n\]
is an isomorphism. To summarize, \m{(T_{ij})_{i\in I,j\in J_i}} (resp. 
\m{(Z_{ij})_{i\in I,j\in J_i}}) is an affine open cover of \m{{\bf B}_{Y,X}} 
(resp. \m{{\bf B}_{Y,X_n}}), and the \m{\rho_{ij}} are the trivializations that 
define the primitive multiple scheme structure on \m{{\bf B}_{Y,X_n}}. The 
restriction of \m{\rho_{ij}},
\[Z_{ij}\cap(X_n\backslash Y)\lra\big(T_{ij}\cap(X\backslash Y)\big)\times{\bf 
Z}_n \ , \]
is \m{\lambda_{i|Z_{ij}\cap(X_n\backslash Y)}}.

Let \m{i,k\in I}, \m{j\in J_i} and \m{l\in J_k}. Let
\[W_{ij,kl} \ = \ D_+^{\BS{\Sigma}_i}({\bf f}_{i,j})\cap 
D_+^{\BS{\Sigma}_k}({\bf f}_{k,l}) \ \subset \ {\bf B}_{Y\cap 
U_i,U_i\times{\bf Z}_n} \ . \]
Then
\[\wt{\lambda_i}^{-1}(W_{ij,kl}) \ = \ \wt{\lambda_k}^{-1}(W_{kl,ij}) \ = \ 
Z_{ij}\cap Z_{kl} \ = \ Z_{ij,kl} \ , \] 
and \ \m{\rho_{ij,kl}^*=\rho_{ij}^{*-1}\rho_{kl}^*} \ is an automorphism of \ 
\m{\ko_{{\bf B}_{Y,X}}(T_{ij,kl})[t]/(t^n)}.

\sepprop

{\em The associated line bundle -- } We have \ 
\m{\rho_{ij,kl}^*(t)=\alpha_{ik}t}, and 
\[\rho_{ij,kl}^*\Big(\frac{t}{{\bf f}_{kl}}\Big) \ = \ \alpha_{ik}\frac{{\bf 
f}_{ij}}{\lambda_{ik}^*\Big({\bf f}_{kl}\Big)}\frac{t}{{\bf f}_{ij}} \ . \]
The coefficient \ \m{\dsp\alpha_{ik}\frac{{\bf 
f}_{ij}}{\lambda_{ik}^*\Big({\bf f}_{kl}\Big)}} \ restricted to \m{T_{ij,kl}} 
is \ \m{\dsp\alpha_{ik|U_{ik}}\frac{{\bf f}_{ij}}{{\bf f}_{kl}}}.
It follows that the line bundle on \m{{\bf B}_{Y,X}} associated to \m{{\bf 
B}_{Y,X_n}} is
\[\ko_{{\bf B}_{Y,X_n}}(\wt{Y}_n)_{|{\bf B}_{Y,X}}\ot\pi_{Y,X}^*(L) \ = \ 
\ko_{{\bf B}_{Y,X}}(\wt{Y})\ot\pi_{Y,X}^*(L) \ . \]

The last statement, about the primitive multiple structure of \m{\wt{Y}_n}, is 
an easy consequence of the previous description of \m{{\bf B}_{Y,X_n}}.
\end{proof}

\sepprop

\begin{subsub}\label{db_sc} The case of double schemes -- \rm We suppose that 
\m{n=2} and keep the notations of the proof of theorem \ref{theo1}. Let 
\m{i,k\in I}. From \ref{DBC}, \m{\lambda_{ik}^*} is an automorphism 
of \m{\ko_X(U_{ik})[t]/(t^2)}, and there exists a derivation \m{D_{ik}} of 
\m{\ko_X(U_{ik})} such that, for every \m{u\in\ko_X(U_{ik})} we have \ 
\m{\lambda_{ik}^*(u)=u+D_{ik}(u)t}.

Let \m{j\in J_i}, \m{l\in J_k}. Then \m{\rho_{ij,kl}^*} is an automorphism of 
\m{\ko_{{\bf B}_{Y,X}}(T_{ij,kl})[z]/(z^2)}, and there is a derivation 
\m{E_{ij,kl}} of \m{\ko_{{\bf B}_{Y,X}}(T_{ij,kl})} such that, for every 
\m{v\in\ko_{{\bf B}_{Y,X}}(T_{ij,kl})}, we have \Nligne 
\m{\rho_{ij,kl}^*(v)=v+E_{ij,kl}(v)z}.

Let \ \m{\Sigma=\{f_{ij}^mf_{kl}^n \ ; \ m,n\geq 0\}}. It is a multiplicative 
system. Let \ \m{R_{ij,kl}=S(U_{ik},V)_\Sigma}, and \m{R_{(ij,kl)}} the ring of 
elements of degree zero in \m{R_{ij,kl}}. We have
\[R_{(ij,kl)} \ = \ \sigg_{m,n\geq 0}I_{mn} \ , \qquad \text{with} \ \
I_{mn}=\Big\{\frac{\sigma_{m+n}}{{\bf f}_{ij}^m{\bf f}_{kl^n}} \ ; \ 
\sigma_{m+n}\in\ki_V(U_{ik})^{m+n}\Big\} \ , \]
and
\[T_{ij,kl} \ = \ \spec(R_{(ij,kl)}) \ . \]
We have then \ \m{E_{ij,kl}={\bf f}_{ij}D_{ik}} (the factor \m{f_{ij}} comes 
from the fact that \ \m{\dsp z=\frac{t}{{\bf f}_{ij}}} \ by step 2).

In this way one obtains a linear map
\[\xi_Y:H^1(X,T_X\ot L)\lra H^1({\bf B}_{Y,X},T_{{\bf 
B}_{Y,X}}\ot\pi_{Y,X}^*(L)\ot\ko_{{\bf B}_{Y,X}}(\wt{Y})) \ , \]
that sends the element of \m{H^1(X,T_X\ot L)} represented by \m{(D_{ij})} to 
the element of \Nligne \m{H^1({\bf B}_{Y,X},T_{{\bf 
B}_{Y,X}}\ot\pi_{Y,X}^*(L)\ot\ko_{{\bf B}_{Y,X}}(\wt{Y}))} represented by 
\m{(E_{ij,kl})}. This map induces the one
\[{\bf DS}(X,L)\lra{\bf DS}({\bf B}_{Y,X},\pi_{Y,X}^*(L)\ot\ko_{{\bf 
B}_{Y,X}}(\wt{Y}))\]
which sends \m{X_2} to \m{{\bf B}_{Y,X_2}}.
\end{subsub}

\sepprop

\begin{subsub}\label{prop7}{\bf Proposition: } Suppose that \m{\codim(Y,X)\geq 
2}. Then \m{\xi_Y} is injective.
\end{subsub}
\begin{proof} Suppose that \m{(E_{ij,kl})} represents 0. Then for every \m{i\in 
I}, \m{j\in J_i}, there exists a derivation \m{F_{ij}} of \m{\ko_{{\bf 
B}_{Y\cap U_i,U_i}}(D_+^{{\bf S}_i}({\bf f}_{ij}))} such that
\[E_{ij,kl} \ = \ F_{ij}-\alpha_{ik}\frac{{\bf f}_{ij}}{{\bf f}_{kl}}F_{kl} \ , 
\]
i.e. \ \m{\dsp D_{ik}=\frac{1}{{\bf f}_{ij}}F_{ij}-\alpha_{ik}\frac{1}{{\bf 
f}_{kl}}F_{kl}}. In particular, if \m{k=i} we have \ \m{\dsp\frac{1}{{\bf 
f}_{ij}}F_{ij}=\frac{1}{{\bf f}_{il}}F_{il}} \ for every \m{j,l\in J_i}. Hence 
there exits a derivation \m{G_i} of  \ \m{\ko_{{\bf B}_{Y\cap 
U_i,U_i}}({\bf B}_{Y\cap U_i,U_i})} \ such that
\ \m{\dsp\frac{1}{{\bf f}_{ij}}F_{ij}=G_{i|D_+^{{\bf S}_i}({\bf f}_{ij})}} 
\ for every \m{j\in J_i}. We have \ \m{U_i\backslash Y\subset{\bf B}_{Y\cap 
U_i,U_i}}, hence \m{G_i} induces a derivation of \m{\ko_X(U_i\backslash Y)}, 
and since \m{\codim(Y,X)\geq 2}, it is also a derivation of \m{\ko_X(U_i)}. We 
have \ \m{D_{ik}=G_i-\alpha_{ik}G_k}. It follows that \m{(D_{ik})} represents 0 
in \m{H^1(X,T_X\ot L)}.
\end{proof}

\end{sub}

\sepsub

\Ssect{Blowing-up of good subschemes}{BLG}

Let \m{Z\subset X_n} be a good subscheme (cf. \ref{GS}) of dimension $d$. Let \ 
\m{Y=Z_{red}} \ and \ \m{\wt{Z}_n=\pi_{Z,X_n}^{-1}(Z)}.

\sepprop

\begin{subsub}\label{theo2}{\bf Theorem: } \m{{\bf B}_{Z,X_n}} is a primitive 
multiple scheme of multiplicity $n$, \m{({\bf B}_{Z,X_n})_{red}={\bf B}_{Y,X}}, 
and its associated line bundle is \ \m{\pi_{Y,X}^*(L)}.

The exceptional divisor \m{\wt{Z}_n} is a primitive multiple scheme of 
multiplicity $n$,\Nligne \m{(\wt{Z}_n)_{red}=\wt{Y}=\pi_{Y,X}^{-1}(Y)}, and 
its associated line bundle is \ \m{\pi_{Y,X}^*(L)_{|\wt{Y}}}.
\end{subsub}
\begin{proof}
It is very similar to that of theorem \ref{theo1}. We will only give 
indications. We use the notations of the proof of this theorem.

From proposition \ref{prop6}, there is an open affine cover \m{(U_i)_{i\in I}} 
of $X$ such that, for every \m{i\in I}, there is a trivialization \ 
\m{\delta_i:U_i^{(n)}\to U_i\times{\bf Z}_n}, \
\m{y_{i,1},\ldots,y_{i,d}\in\ki_{Y,X}(U_i)} \ that generate \m{\ki_{Y\cap 
U_i,U_i}} and \m{\ki_{\delta_i(Z\cap U_i^{(n)}),U_i\times{\bf Z}_n}}.

Let \m{U\subset X} be a smooth open subset, and \m{V\subset U} a closed smooth
connected subvariety of codimension $d$, such that there exist generators 
\m{y_1,\ldots,y_d} of \m{\ki_{V,U}}. Let \ \m{U_n=U\times{\bf Z}_n}. We have \ 
\m{U_n=\spec(\ko_X(U)[t]/(t^n))}. Let \m{Z\subset U_n} be the subscheme 
corresponding to the ideal \ \m{(y_1,\ldots,y_d)\subset\ko_X(U)[t]/(t^n)}.

We have \ \m{{\bf B}_{Z,U_n}=\Proj(\B_n)}, where \ \m{\B_n=\B_n(Z,U)} \ is 
the graduate ring
\[\B_n \ = \ \bigoplus_{p\geq 0}\ki_{Z,U_n}(U)^p \ , \]
and
\[\ki_{Z,U_n}(U)^p \ = \ (y_1,\ldots,y_d)^p \ . \]
As in step 2 of the proof of theorem \ref{theo1},  if \ \m{f\in\ki_{Z,U_n}(U)}, 
we can define the affine open subset \m{D^{\B_n}_+(f)} of \m{{\bf B}_{Z,U_n}}. 
Let \m{\B_{n,f}} be the localized ring, and \ \m{\B_{n,(f)}\subset\B_{n,f}} \ 
the subring of elements of degree 0. Then \ 
\m{D^{\B_n}_+(f)\simeq\spec(\B_{n,(f)})}. We have 
\[\B_{n,(f)} \ = \ \sum_{p\geq 0}\B_{n,(f),p} \ , \quad \text{with} \ \ \
\B_{n,(f),p} \ = \ \Big\{\frac{\theta_p}{f^p}\ ;\ 
\theta_p\in\ki_{Z,U_n}(U)^p\Big\} \ . \]
It follows that \ \m{\B_{n,(f)}=S_{(f)}[t]/(t^n)}. We obtain an isomorphism \ 
\m{D^{\B_n}_+(f)\simeq D^S_+(f)\times{\bf Z}_n} \ inducing the identity on 
\m{D^S_+(f)}. 

Using the open cover \m{(D^{\B_n}_+(y_i))_{i\in I}}, it is easy to see that we 
have a canonical isomorphism \ \m{{\bf B}_{Z,U_n}\simeq{\bf B}_{Y,U}\times{\bf 
Z}_n}. It is now straightforward to follow the same way as in steps 5,6 of the 
proof of theorem \ref{theo1}, and finish the proof of theorem \ref{theo2}.
\end{proof}

\end{sub}

\sepsub

\Ssect{Successive blowing-ups}{SBU}

Let \m{Z\subset X_n} be a good subscheme (cf. \ref{GS}). Let \m{Y=Z_{red}}, and 
\m{\wt{Y}=\pi_{Y,X}^{-1}(Y)\subset{\bf B}_{Y,X}}. Then \m{\wt{Y}} is a Cartier 
divisor of \m{{\bf B}_{Y,X}}. Hence \ \m{{\bf B}_{\wt{Y},{\bf B}_{Y,X}}={\bf 
B}_{Y,X}}.

\sepprop

The following result is an easy consequence of the descriptions of the 
blowing-ups in theorems \ref{theo1} and \ref{theo2}.

\sepprop

\begin{subsub}\label{prop11}{\bf Proposition: } The blowing-up of \m{{\bf 
B}_{Z,X_n}} along $\wt{Y}$ is canonically isomorphic to \m{{\bf B}_{Y,X_n}}, 
i.e. \ \m{{\bf B}_{\wt{Y},{\bf B}_{Z,X_n}}\simeq{\bf B}_{Y,X_n}}.
\end{subsub}

\end{sub}

\sepsec

\section{Blowing-up of good zero-dimensional subschemes of primitive 
double schemes}\label{pt}

We use the notations of \ref{BUSCH}. We suppose that \ \m{d=\dim(X)\geq 2}, 
that $Y$ is a single point: \m{Y=\{P\}}, and that \m{n=2}. There exists a 
finite open affine cover \m{(U_i)_{i\in I}} of $X$ such that
\begin{enumerate}
\item[--] There is a unique \m{i_0\in I} such that \m{P\in U_{i_0}}, and \ 
\m{\cup_{i\in I,i\not=i_0}U_i=X\backslash\{P\}}.
\item[--] $\ki_{\{P\},U_{i_0}}$ is generated by $d$ elements \m{y_1,\ldots,y_d},
and \m{\Omega_{X|U_{i_0}}} by \m{dy_1,\ldots,dy_d}.
\item[--] $L$ is defined by a cocycle \m{(\alpha_{ij})}, with \ 
\m{\alpha_{ij}\in\ko_X(U_{ij})^*}.
\item[--] For every \m{i\in I}, there is an isomorphism \ 
\m{\delta_i:U_i^{(2)}\to U_i\times{\bf Z}_2} (where \m{U_i^{(2)}} is the open 
subscheme of \m{X_2} corresponding to \m{U_i}), inducing the identity on 
\m{U_i}. There is a derivation \m{D_{ij}} of \m{\ko_X(U_{ij})} such that, for 
every \m{u\in\ko_X(U_{ij})}, \m{\delta_{ij}^*(u)=u+D_{ij}(u)t}. We have \ 
\m{\delta_{ij}^*(t)=\alpha_{ij}t}.
\end{enumerate}
We consider the open affine cover of \m{{\bf B}_{P,X}} which consists of the 
\m{U_i}, \m{i\in I}, and the \m{D_+^{\B_2}(y_k)\subset{\bf B}_{P,U_{i_0}}}, 
\m{1\leq 
k\leq d}. Let \ \m{J=I\cup\{i_1,\ldots,i_d\}}, \m{V_i=U_i} \ if \m{i\in I} and 
\ \m{V_{i_k}=D_+^{\B_2}(y_k)}. The line bundle \m{\pi_{P,X}^*(L)} on \m{{\bf 
B}_{P,X}} is defined by the cocycle \m{(\beta_{ij})_{i,j\in J}}, 
\m{\beta_{ij}\in\ko_{{\bf B}_{P,X}}(V_{ij})^*}, with \ 
\m{\beta_{ij}=\alpha_{ij}} \ if \m{i,j\in I}, \m{\beta_{ki_m}=\alpha_{ki_0}}, 
\m{\beta_{i_mi_n}=1} \ if \m{k\in I}, \m{1\leq m,n\leq d}.

We will study the blowing-ups of \m{X_2} along good subschemes $Z$ such that 
\m{Z_{red}=\{P\}}. From theorem \ref{theo2}, \m{{\bf B}_{Z,X_2}} is a primitive 
multiple scheme of multiplicity $2$, \m{({\bf B}_{Z,X_2})_{red}={\bf B}_{P,X}}, 
and its associated line bundle is \m{\pi_{P,X}^*(L)}. The exceptional divisor 
\m{\wt{P}=\pi_{P,X}^{-1}(P)} is a projective space of dimension \m{d-1}. Let 
\m{\wt{P}_2=\pi^{-1}_{Z,X_2}(P)}. It is the trivial primitive multiple scheme 
of multiplicity $2$, such that \m{(\wt{P}_2)_{red}=\wt{P}}, and with associated 
line bundle \ \m{\ko_{\wt{P}}}.

\sepsub

\Ssect{Parametrization of good zero-dimensional subschemes}{par_go}

We will prove

\sepprop

\begin{subsub}\label{prop10}{\bf Proposition: }Let \m{\BS{\Gamma}_P} be the set 
of good subschemes \m{Z\subset X_2} such that \ \m{Z_{red}=\{P\}}. Then 
\m{\BS{\Gamma}_P} has a natural structure of affine space, with associated 
vector space \m{(T_X\ot L)_P}.
\end{subsub}

\sepprop

Let \m{\mm_P} be the maximal ideal of \m{\ko_{X,P}}. Let \ 
\m{a_1,\ldots,a_d\in\ko_{X,P}}, \Nligne
\m{I(a_1,\ldots,a_d)=(y_1+ta_1,\ldots,y_d+ta_d)}, and \m{Z_0(a_1,\ldots,a_d)} 
the corresponding subscheme of \m{U_{i_0}\times{\bf Z}_2}. The good 
subschemes \ \m{Z\subset U_{i_0}\times{\bf Z}_2} \ such that \ 
\m{Z_{red}=\{P\}} \ are precisely the \m{Z_0(a_1,\ldots,a_d)}.

 Let \ \m{Z(a_1,\ldots,a_d)=\delta_{i_0}^{-1}(Z_0(a_1,\ldots,a_d))\subset X_2}. 
It is the subscheme defined by the ideal  sheaf \ 
\m{\delta_{i_0}^*(I(a_1,\ldots,a_d))}. The good subschemes \ \m{Z\subset X_2} \ 
such that \ \m{Z_{red}=\{P\}} \ are the \m{Z(a_1,\ldots,a_d)}.

\sepprop

\begin{subsub}\label{lem7}{\bf Lemma: } Let 
\m{a_1,a'_1,\ldots,a_d,a'_d\in\ko_{X,P}}. Then \ \m{I(a_1,\ldots,a_d)\simeq 
I(a'_1,\ldots,a'_d)} \ if and only if \ \m{a_i-a'_i\in\mm_P} \ for \m{1\leq 
i\leq d}.
\end{subsub}
\begin{proof}
Suppose that \ \m{I(a_1,\ldots,a_d)\simeq I(a'_1,\ldots,a'_d)}. Then there 
exists an invertible \m{(d\times d)}-matrix \m{M=(\theta_{ij})} of elements of 
\m{\ko_{U_{i_0}\times{\bf Z}_2,P}} such that 
\[\begin{pmatrix}y_1+ta'_1\\.\\.\\y_d+ta'_d\end{pmatrix}=
M.\begin{pmatrix}y_1+ta_1\\.\\.\\y_d+ta_d\end{pmatrix} \ .
\]
We can write \ \m{\theta_{ij}=\lambda_{ij}+t\mu_{ij}}, with 
\m{\lambda_{ij},\mu_{ij}\in\ko_{X,P}}. We have then
\[y_i \ = \ \sigg_{j=1}^d\lambda_{ij}y_j \ , \qquad
a'_i \ = \ \sigg_{j=1}^d\lambda_{ij}a_j+\sigg_{j=1}^d\mu_{ij}y_j \ . \]
The Koszul relations imply that we have \ \m{\lambda_{ii}=1}, and that \ 
\m{\lambda_{ij}\in\mm_P} \ if \m{i\not=j}. It follows that \ 
\m{a_i-a'_i\in\mm_P}.

Conversely, if \m{a_i-a'_i\in\mm_P} \ for \m{1\leq i\leq d}, we can write \ 
\m{a'_i=a_i+\sigg_{j=1}^d\rho_{ij}y_j}, \m{\rho_{ij}\in\ko_{X,P}}. We have then
\[y_i+ta'_i \ = \ y_i+ta_i+\sigg_{j=1}^dt\rho_{ij}(y_j+ta_j) \ , \]
hence \ \m{I(a_1,\ldots,a_d)\simeq I(a'_1,\ldots,a'_d)}.
\end{proof}

\sepprop

It follows from lemma \ref{lem7} that the map
\[\xymatrix@R=1pt{\C^d\ar[r] & \BS{\Gamma}_P \\ 
(\zeta_1,\ldots,\zeta_d)\fmaps[r] & Z(\zeta_1,\ldots,\zeta_d)}\]
is a bijection.

Let \m{U\subset X} be an open affine neighborhood of $P$, and $D$ a derivation 
of \m{\ko_X(U)}, that we can view as a section of \m{T_{X|U}}. Let \ \m{D_P\in 
T_{X,P}} \ its value at $P$. Let \m{a_1,\ldots,a_d\in\ko_{X,P}}. Then \ 
\m{\dsp D(a_1,\ldots,a_d)=\sigg_{i=1}^da_i\frac{\partial}{\partial y_i}} \ is a 
derivation of \m{\ko_X(U)}, where $U$ is a suitable neighborhood of $P$. We have
\[D(a_1,\ldots,a_d)_P \ = \ \sigg_{i=1}^da_i(P)\Big(\frac{\partial}{\partial 
y_i}\Big)_P \ . \]
According to lemma \ref{lem7}, \m{Z(a_1,\ldots,a_d)=Z(a'_1,\ldots,a'_d)} \ if 
and only if \Nligne \m{D(a_1,\ldots,a_d)_P=D(a'_1,\ldots,a'_d)_P}.

If \m{Z\in\BS{\Gamma}_P}, and \ \m{Z=Z(a_1,\ldots,a_d)}, let
\[\Delta_P(Z,(y_i)_{1\leq i\leq d},\delta_{i_0}) \ = \ D(a_1,\ldots,a_d)_P
\ot\delta_{i_0}^*(t)(P) \ \in \ (T_X\ot L)_P \ . \]

\sepprop

\begin{subsub}\label{lem8}{\bf Lemma: } \m{\Delta_P(Z,(y_i)_{1\leq i\leq 
d},\delta_{i_0})} does not depend on \m{(y_i)_{1\leq i\leq d}}.
\end{subsub}
\begin{proof}
Let \m{z_1,\ldots,z_d} be generators of \m{\mm_P}. We have \ 
\m{(y_1,\ldots,y_d)=(z_1,\ldots,z_d)M}, where  \ \m{M=(\mu_{ij})_{1\leq i,j\leq 
d}} \ is an invertible \m{(d\times d)}-matrix of elements of \m{\ko_{X,P}}. We 
have
\[\deriv{y_i}{z_k} \ = \ \mu_{ki}+\sigg_{j=1}^d\deriv{\mu_{ji}}{z_k}z_j \ , \]
hence
\[\Big(\Big(\deriv{}{z_1}\Big)_P,\ldots,\Big(\deriv{}{z_d}\Big)_P\Big) \ = \ 
{\Big(\Big(\deriv{}{y_1}\Big)_P,\ldots,\Big(\deriv{}{y_d}\Big)_P\Big)}.^tM(P) 
\ . \]
Let \ \m{(b_1,\ldots,b_d)=(a_1,\ldots,a_d)M(P)^{-1}}. Then
\[(y_1+a_1t,\ldots,y_d+a_dt) \ = \ (z_1+b_1t,\ldots,z_d+b_dt)\]
in \m{\ko_{X,P}[t]/(t^2)}. We have
\begin{eqnarray*}
\Delta(Z,(z_i)_{1\leq i\leq d},\delta_{i_0}) & = & 
\Big(\Big(\deriv{}{z_1}\Big)_P,\ldots,\Big(\deriv{}{z_d}\Big)_P\Big)
\begin{pmatrix}
b_1 \\ . \\ . \\ b_d
\end{pmatrix}\ot\delta_{i_0}^*(t)(P)\\
& = & \Big(\Big(\deriv{}{y_1}\Big)_P,\ldots,\Big(\deriv{}{y_d}\Big)_P\Big).
{}^tM(P).^tM(P)^{-1}\begin{pmatrix}
a_1 \\ . \\ . \\ a_d
\end{pmatrix}\ot\delta_{i_0}^*(t)(P)\\
& = & \Delta(Z,(y_i)_{1\leq i\leq d},\delta_{i_0}) \ .
\end{eqnarray*}
\end{proof}

\sepprop

So we will omit the parameter \m{(y_i)_{1\leq i\leq d}}, and note \ 
\m{\Delta(Z,(y_i)_{1\leq i\leq d},\delta_{i_0})=\Delta(Z,\delta_{i_0})}. The 
following lemma ends the proof of proposition \ref{prop10}.

\sepprop

\begin{subsub}\label{lem10}{\bf Lemma} Let \m{Z,Z'\in\BS{\Gamma}_P} and \ 
\m{\phi:U_{i_0}^{(2)}\to U_{i_0}\times{\bf Z}_2} \ an isomorphism leaving 
\m{U_{i_0}} invariant. Then we have
\[\Delta_P(Z',\delta_{i_0})-\Delta_P(Z,\delta_{i_0}) \ = \ 
\Delta_P(Z',\phi)-\Delta_P(Z,\phi) \ . \]
\end{subsub}
\begin{proof}
Let \ \m{\theta=\phi\circ\delta_{i_0}^{-1}}, which is an automorphism of 
\m{U_{i_0}\times{\bf Z}_2}. We have
\[\theta^{*-1}(y_i)=y_1+\epsilon_it \ \ \text{for} \ \ 1\leq i\leq d \ , 
\qquad \theta^{*-1}(t) \ = \ \alpha t \ , \]
with \ \m{\epsilon_i\in\ko_X(U_{i_0})}, \m{\alpha\in\ko_X(U_{i_0})^*}. Suppose 
that \ \m{Z=Z(a_1,\ldots,a_d)}. Then
\[Z \ = \ \delta_{i_0}^{-1}(Z(a_1,\ldots,a_d)) \ = \ 
\phi^{-1}\big(\theta(Z(a_1,\ldots,a_d))\big) \ , \]
and \m{\theta(Z(a_1,\ldots,a_d))} is the subscheme whose ideal is generated by 
\m{\theta^{*-1}(y_i+\epsilon_it)}, \m{1\leq i\leq t}. We have
\[\theta^{*-1}(y_i+a_it) \ = \ y_i+(\epsilon_i+\alpha a_i)t \ , \quad
\text{hence} \ \ \delta_{i_0}^*(y_i+a_it) \ = \ \phi^*(y_i+(\epsilon_i+\alpha 
a_i)t) \ , \]
and \ \m{\dsp\phi^*(t)=\delta_{i_0}^*(\theta^*(t))=\frac{1}{\alpha}
\delta_{i_0}^*(t)}. It follows that
\[\Delta_P(Z,\phi) \ = \ \Big(\sigg_{i=1}^d(\epsilon_i(P)+\alpha(P)a_i(P))
\Big(\deriv{}{y_i}\Big)_P\Big)\ot
\Big(\frac{1}{\alpha(P)}\delta_{i_0}^*(t)(P)\Big) \ . \]
We have a similar formula for \ \m{Z'=Z(a'_1,\ldots,a'_d)}, and the result 
follows easily.
\end{proof}

From proposition \ref{prop10} and the preceding lemmas, we have \m{Z=Z'} if and 
only if \Nligne \m{\Delta_P(Z,\delta_{i_0})=\Delta_P(Z',\delta_{i_0})}.

\sepprop

\begin{subsub}\label{cas_triv} The case of a trivial primitive double scheme 
-- \rm If \m{X_2} is trivial, we can choose the trivializations such that \ 
\m{D_{ij}=0} for every \m{i,j\in I}. We obtain a canonical identification \ 
\m{\BS{\Gamma}_P=(T_X\ot L)_P}.

In this case \m{X_2} is the second infinitesimal neighborhood of $X$ in the 
line bundle \m{L^*} (viewed as a smooth variety, and \m{X\subset X_2} defined 
by the zero section of \m{L^*}). We have a canonical projection \ \m{p:X_2\to 
X}.

Let \m{Z_0\subset X_2} be the good subscheme corresponding to \ \m{0\in(T_X\ot 
L)_P}, it is \m{p^{-1}(P)}. Then \m{{\bf B}_{Z_0,X_2}} is the trivial primitive 
scheme with associated smooth variety \m{{\bf B}_{P,X}} and associated line 
bundle \m{\pi_{P,X}^*(L)}. 
\end{subsub}

\end{sub}

\sepsub

\Ssect{Blowing-ups}{BL}

\begin{subsub}\label{lem11}{\bf Lemma: } Let $E$ be a derivation of 
\m{U_{i_0}\backslash\{P\}}. Then $U$ can be extended uniquely to a derivation 
of \m{U_{i_0}}. It can be extended to a derivation of \m{\pi_{P,{\bf 
B}_{Z,X_2}}^{-1}(U_{i_0})} if and only if \m{E_P=0}.
\end{subsub}
\begin{proof}
Since \ \m{\dim(X)\geq 2}, $E$ can be extended to \m{U_{i_0}}. It can be 
extended to \m{\pi_{P,{\bf B}_{Z,X_2}}^{-1}(U_{i_0})} if and only, for \m{1\leq 
i\leq d}, if it be extended to \m{D_+^{\B_n}(y_i)}. For every 
\m{\sigma\in\ki_{P,U_{i_0}}}, we have
\[E\Big(\frac{\sigma}{y_i}\Big) \ = \ \frac{E(\sigma)y_i-E(y_i)\sigma}{y_i^2} \ 
. \]
Hence $E$ can extended to \m{D_+^{\B_n}(y_i)} if and only if \ 
\m{E(\sigma)y_i-E(y_i)\sigma\in\ki_{P,U_{i_0}}^2}, for every 
\m{\sigma\in\ki_{P,U_{i_0}}}. It is easy to see that this happens if and only 
if \m{E_P=0}.
\end{proof}

\sepprop

Let \ \m{\text{ev}_L:\ko_X\ot H^0(X,T_X\ot L)\to T_X\ot L} \ be the evaluation 
morphism, and \Nligne \m{H_{L,P}=\imm(\text{ev}_{L,P})\subset(T_X\ot L)_P}. Let 
\ \m{\Psi_{L,P}:T_{X_P}\ot L\to (T_{X,P}\ot L)/H_{L,P}} \ be the quotient map.

Let \ \m{Z,Z'\in\BS{\Gamma}_P}. Then \m{{\bf B}_{Z,X_2}} and \m{{\bf 
B}_{Z',X_2}} are primitive double schemes with associated smooth variety 
\m{{\bf B}_{P,X}}. We say that they are {\em isomorphic} if there exists an 
isomorphism \ \m{{\bf B}_{Z,X_2}\simeq{\bf B}_{Z',X_2}} \ inducing the identity 
on \m{{\bf B}_{P,X}}.

\sepprop

\begin{subsub}\label{theo3}{\bf Theorem: } If \m{X_2} is trivial, we choose the 
trivializations such that \ \m{D_{ij}=0} for every \m{i,j\in I}. Let \ 
\m{Z,Z'\in\BS{\Gamma}_P}. 

{\bf 1 -- } If \m{X_2} is not trivial, \m{{\bf B}_{Z,X_2}} and \m{{\bf  
B}_{Z',X_2}} are isomorphic if and only if \Nligne 
\m{\Delta_P(Z,\delta_{i_0})-\Delta_P(Z',\delta_{i_0})\in H_{L,P}}.

{\bf 2 -- } If \m{X_2} is trivial, \m{{\bf B}_{Z,X_2}} and \m{{\bf  
B}_{Z',X_2}} are isomorphic if and only if \Nligne
\m{\C.\Psi_{L,P}(\Delta_P(Z',\delta_{i_0}))=
\C.\Psi_{L,P}(\Delta_P(Z,\delta_{i_0}))}.
\end{subsub}
\begin{proof} According to \ref{BLG}, the double structures on 
\m{{\bf B}_{Z,X_2}} and \m{{\bf B}_{Z',X_2}} are defined by families 
\m{(E_{ij})_{i,j\in J}}, \m{(E'_{ij})_{i,j\in J}}, where \m{E_{ij}} and 
\m{E'_{ij}} are derivations of \m{\ko_X(V_{ij})}. These families define the 
elements \m{h_Z,h_{Z'}} of \m{H^1(X,T_{{\bf B}_{P,X}}\ot L)} associated to 
\m{{\bf B}_{Z,X_2}} and \m{{\bf B}_{Z',X_2}} respectively. 

Suppose that \ \m{Z=Z(a_0,\ldots,a_d)}, \m{Z'=Z(a'_1,\ldots,a'_d)}. Let \ 
\m{D=D(a_0,\ldots,a_d)},\Nligne \m{D'=D(a'_0,\ldots,a'_d)}. Suppose that 
\m{Z\not=Z'}. Then \ \m{D_P\not=D'_P}. We can assume that \ \m{D'_P\not=0}. By 
adding to $D$ a suitable derivation \m{F} such that \m{F_P=0}, we can assume 
that, for every \ \m{\lambda\in\C}, we have \ \m{D\not=\lambda D'}.

Let $\theta$ be the automorphism of \ \m{U_{i_0}\times{\bf Z}_2} \ defined by 
the derivation $D$, i.e. \ \m{\theta^*(u)=u+D(u)t} \ for every 
\m{u\in\ko_X(U_{i_0})}, and \m{\theta^*(t)=t}. We have \ 
\m{\theta^*(y_i)=y_i+a_it} \ for \m{1\leq i\leq d}. Hence 
\m{\theta\circ\delta_{i_0}(Z)} is the subscheme defined by the ideal 
\m{(y_1,\ldots,y_d)}. If instead of \m{\delta_{i_0}} we 
take the trivialization \m{\theta\circ\delta_{i_0}}, and do a similar change 
for \m{Z'}, we obtain by \ref{ef_ch} and the proof of theorem \ref{theo2} that
\[E_{ij}=E'_{ij}=D_{ij} \quad \text{if} \ \ i,j\in I\backslash\{i_0\} \ , 
\qquad \ E_{i_mi_n}=E'_{i_mi_n}=0 \quad \text{if} \ \ 1\leq m,n\leq d \ , \]
\[E_{ki_m}=D_{ki_0}+\alpha_{ki_0}D \ , \quad E'_{ki_m}=D_{ki_0}+\alpha_{ki_0}D' 
\qquad \text{if} \ \ k\in 
I\backslash\{i_0\}, \ 1\leq m\leq d \ . \]
From \ref{DBC}, \m{{\bf B}_{Z,X_2}} and \m{{\bf B}_{Z',X_2}} are isomorphic if 
and only if there exists \m{\lambda\in\C^*} such that \ \m{h_Z=\lambda h_{Z'}}, 
i.e. if and only if there exist, for every \m{j\in J}, a derivation \m{G_j} of 
\m{\ko_{{\bf B}_{P,X}}(V_j)} such that \ \m{E_{ij}=\lambda 
E'_{ij}+G_i-\beta_{ij}G_j} \ on \m{V_{ij}}, for all \m{i,j\in J}. These 
relations are equivalent to \ \m{G_{i_m}=G_{i_n}} \ on \m{V_{i_mi_n}}, for 
every \m{i,j\in I\backslash\{i_0\}}
\begin{equation}\label{equ2}(1-\lambda)D_{ij} \ = \ G_i-\alpha_{ij}G_j 
\qquad \text{on} \ \ U_{ij} \ ,
\end{equation}
and for \m{k\in I\backslash\{i_0\}}, \m{1\leq m\leq d},
\begin{equation}\label{equ3} (1-\lambda)D_{ki_0} \ = \ \alpha_{ki_0}(\lambda 
D'-D)+G_k-\alpha_{ki_0}G_{i_m} \qquad \text{on} \ \ V_{ki_m} \ .
\end{equation}
The derivations \m{G_{i_m}} do not depend on $m$, and can be extended to a 
derivation $\Theta$ on \m{U_{i_0}}, as well as on \m{\pi_{P,{\bf 
B}_{Z,X_2}}^{-1}(U_{i_0})}, and by lemma \ref{lem11} we have \ \m{\Theta_P=0}. 

Suppose that \ \m{\lambda=1}. Then from $(\ref{equ2})$, and the fact that \ 
\m{\cup_{i\in I,i\not=i_0}U_i=X\backslash\{P\}}, there exists \ \m{\Phi\in 
H^0(X,T_X\ot L)} such that \ \m{G_i\ot\alpha_i^{-1}=\Phi_{|U_i\ot L}} \ for 
every \m{i\in I\backslash\{i_0\}}. From $(\ref{equ3})$, we have 
\[D'-D-\Theta+\Phi\ot\alpha_{i_0} \ = 0 \ , \]
which implies that \ \m{\Delta_P(Z,\delta_{i_0})-\Delta_P(Z',\delta_{i_0})\in 
H_{L,P}}.

Suppose that \m{\lambda\not=1}. Let \ \m{\dsp G_{i_0}=\Theta+D-\lambda D'}, and 
for every \m{i\in I}, \m{\dsp F_i=\frac{1}{1-\lambda}G_i}. Then we have, for 
every \m{i,j\in I\backslash\{i_0\}}
\[D_{ij} \ = \ F_i-\alpha_{ij}F_j \qquad \text{on} \ \ U_{ij} \ , \]
and, for every \m{k\in I\backslash\{i_0\}}, \m{1\leq m\leq d}
\[D_{ki_0} \ = \ F_k-\alpha_{ki_0}F_{i_0} \qquad \text{on} \ \ V_{ki_m} \ . \]
Hence \ \m{D_{ki_0}=F_k-\alpha_{ki_0}F_{i_0}} \ on \m{U_{ki_0}}. It follows 
that \m{X_2} is trivial. But in this case we have assumed that \m{D_{ij}=0} for 
every \m{i,j\in I}. Hence there exists  \ \m{\Phi\in 
H^0(X,T_X\ot L)} such that \ \m{G_i\ot\alpha_i^{-1}=\Phi_{|U_i\ot L}} \ for 
every \m{i\in I}, and we have
\[\lambda D'-D-\Theta+\Phi\ot\alpha_{i_0} \ = 0 \ , \]
which implies that \ \m{\C.\Psi_{L,P}(\Delta_P(Z',\delta_{i_0}))=
\C.\Psi_{L,P}(\Delta_P(Z,\delta_{i_0}))}.

The converse is immediate.
\end{proof}

\end{sub}

%\sepsec
\newpage

\section{Blowing-up of an hypersurface}\label{BUSHYP}

\Ssect{Primitive multiple rings}{PMR}

Let $R$ be an entire noetherian ring and $K$ its fraction field.

Let $n$ be a positive integer and \ \m{R[n]=R[t]/(t^n)}, which we call the {\em 
primitive multiple ring of multiplicity} $n$ associated to $R$.

An element $u$ of $R[n]$ can be written in an unique way as
\[u \ = \ \sigg_{i=0}^{n-1}u_it^i \ , \]
with \m{u_i\in R} for \ \m{0\leq i\leq n}. Then $u$ is a zero divisor if and 
only if \m{u_0=0}, i.e. if $u$ is a multiple of $t$, and it is invertible if 
and only if \m{u_0} is invertible.

Let \m{S_n} be the multiplicative system of non zero divisors of \m{R[n]}, i.e. 
\ \m{S_n=\{u\in R[n];u_0\not=0\}}. It is easy to see that \ 
\m{S_n^{-1}R[n]=K[t]/(t^n)=K[n]} (the primitive multiple ring of multiplicity 
$n$ associated to $K$).

Suppose that \m{n>1}. Let
\[\xymatrix@R=1pt{\sigma^R_n:R[n]\ar[r] & R[n-1]\\ \ \ u\fmaps[r] & 
\sigg_{i=0}^{n-2}u_it^i } \ , \qquad \qquad
\xymatrix@R=1pt{\rho^R_n:R[n]\ar[r] & R\\ \ \ u\fmaps[r] & u_0 \ .} 
\] 
which are morphisms of rings. 

Let \m{\Gamma_n(R)} be the set of morphisms \m{R\to 
R[n]} such that \ \m{\rho^R_n\circ\phi=I_R}, i.e. such that \ \m{\phi(a)_0=a} \ 
for every \m{a\in R}, and \m{\Sigma_n(R)} the set of endomorphisms $\Theta$ of 
\m{R[n]} such that \ \m{\rho_n^R\circ\Theta=\rho_n^R}.

Let \m{\epsilon\in R[n-1]}, \m{\phi\in\Gamma_n(R)}, and
\[\xymatrix@R=1pt{\Theta^R_n(\phi,\epsilon):R[n]\ar[r] & R[n]\\ \quad
u=\sigg_{i=0}^{n-1}u_it^i\fmaps[r] & \sigg_{i=0}^{n-1}\phi(u_i)\epsilon^it^i \ 
. }\]

\sepprop

\begin{subsub}\label{prop1}{\bf Proposition: } {\bf 1 -- } 
\m{\Theta^R_n(\phi,\epsilon)} is an endomorphism.

{\bf 2 -- } For every \m{\Theta\in\Sigma_n(R)}, there exist a unique 
\m{\phi\in\Gamma_n(R)} and a unique \m{\epsilon\in R[n-1]} such that \ 
\m{\Theta=\Theta^R_n(\phi,\epsilon)}.

{\bf 3 -- } \m{\Theta^R_n(\phi,\epsilon)} is injective if and only \ 
\m{\epsilon_0\not=0}.

{\bf 4 -- } \m{\Theta^R_n(\phi,\epsilon)} is surjective if and only it is an 
isomorphism, if and only if \m{\epsilon_0} is invertible.
\end{subsub}
\begin{proof}
{\bf 1} is immediate. 

Now we prove {\bf 2}. We have \ \m{\Theta(t)^n=\Theta(t^n)=0}, hence 
\m{\Theta(t)} is a multiple ot $t$: \m{\Theta(t)=\epsilon t}, for a unique 
\m{\epsilon\in R[n-1]}. If we take \ \m{\phi=\Theta_{|R}}, it is immediate that 
\ \m{\Theta=\Theta^R_n(\phi,\epsilon)} \ and that $\phi$ is unique.

Suppose that \m{\epsilon_0=0}, i.e. $\epsilon$ is a multiple of $t$. Then \ 
\m{\Theta^R_n(\phi,\epsilon)(t^{n-1})=0}, hence \m{\Theta^R_n(\phi,\epsilon)} 
is not injective. Conversely, suppose that \m{\epsilon_0\not=0}. We will prove 
that \m{\Theta^R_n(\phi,\epsilon)} is injective by induction on $n$. It is 
obvious for \m{n=1}. Suppose that it is true for \m{n-1\geq 1}. Let \ 
\m{\Theta=\Theta^R_n(\phi,\epsilon)}, 
\m{\Theta'=\Theta_{n-1}^R(\sigma_n^R\circ\phi,
\sigma_n^R(\epsilon))\in\Sigma_{n-1}(R)}. Then \m{\Theta'} is injective, and we 
have \ \m{\Theta'=\rho_n^R\circ\Theta_{|R[n-1]}}. Let \m{u\in R[n]}, 
\m{u=v+u_{n-1}t^{n-1}}, \m{v\in R[n-1]}. Suppose that \ \m{\Theta(u)=0}. Then \ 
\m{\Theta(v)=-u_{n-1}\epsilon^{n-1}t^{n-1}}, hence \ \m{\Theta'(v)=0}, and \ 
\m{v=0}. So we have \ \m{\Theta(u)=u_{n-1}\epsilon^{n-1}t^{n-1}=0}. Since \ 
\m{\epsilon_0\not=0}, \m{u_{n-1}=0} \ and \ \m{u=0}. Hence $\Theta$ is 
injective. This proves {\bf 3}. The proof of {\bf 4} is similar.
\end{proof}

\sepprop

Let \m{\phi\in\Gamma_n(R)} (resp. \m{\Theta\in\Sigma_n(R)}). Then $\phi$ 
naturally induces \m{\ov{\phi}\in\Gamma_n(K)} (resp. 
\m{\ov{\Theta}\in\Sigma_n(K)}) such that \ \m{\ov{\phi}_{|R}=\phi} (resp. 
\m{\ov{\Theta}(R[n])\subset R[n]} and \ \m{\ov{\Theta}_{|R[n]}=\Theta}). Let 
\m{\epsilon\in R[n]}. We have
\[\Theta^K_n(\ov{\phi},\epsilon) \ = \ \ov{\Theta^R_n(\phi,\epsilon)} \ , \]
and it is an isomorphism if and only if \ \m{\epsilon_0\not=0}.

For every \ \m{\lambda=\sigg_{i=0}^{n-1}\lambda_it^i\in R[n]} \ and 
\m{\alpha\in R[n]}, let \ \m{\lambda[\alpha t]=\sigg_{i=0}^{n-1}\lambda_i
\alpha^it^i}. If \m{\phi\in\Gamma_n(R)}, let
\[\xymatrix@R=1pt{\phi[\alpha t]:R\ar[r] & R[n]\\ \qquad \ \ 
a\fmaps[r] & \phi(a)[\alpha t] \ , }\]
which also belongs to \m{\Gamma_n(R)}.

Let \m{y\in R}, \m{y\not=0}. Let \ \m{i_R:R\to R[n]} \ be the inclusion, and
\[\chi_y^R \ = \ \Theta(i_R,y) \ . \]
We have \ \m{\ov{\chi_y^R}=\chi_y^K}.

\sepprop

\begin{subsub}\label{prop2}{\bf Proposition: } Let \m{x\in R} such that 
\m{x\not=0}, and \m{\alpha\in R[n]}. Let \
\m{\Theta=\Theta^R_n(\phi,\epsilon)\in\Sigma_n(R)} and
\[\Psi \ = \ \chi_{\alpha x}^K\circ\ov{\Theta}
\circ\chi^K_{1/x} \ \in \ \Sigma_n(K) \ . \]
Let \ \m{\phi(x)=x+\mu t}, with \m{\mu\in R[n-1]}. Then we have \ 
\m{\Psi=\ov{\Theta(\phi',\epsilon')}}, with
\begin{equation}\label{equ1}
\phi' \ = \ \phi[\alpha xt] \ , \quad \epsilon' \ = \ 
\frac{\alpha\epsilon[\alpha xt]}{1+\alpha\mu[\alpha xt]t} \ . \end{equation}
Moreover, if $\Theta$ is an isomorphism and $\alpha$ is invertible, then 
\m{\Theta(\phi',\epsilon')} is also an isomorphism.
\end{subsub}
\begin{proof} Let \m{a\in R}. We have
\[
\Psi(a)=\chi_{\alpha x}^K\circ\ov{\Theta}(a)=
\chi_{\alpha x}^K(\phi(a))=\phi[a][\alpha xt] \ \in \ R[n] \ , \]
and
\[\Psi(t) \ = \ \chi_{\alpha x}^K\circ\ov{\Theta}\Big(\frac{1}{x}t\Big) \ = \
\chi_{\alpha x}^K\Big(\frac{\epsilon}{\phi(x)}t\Big) \ =
\frac{\epsilon[\alpha xt]\alpha x}{\phi(x)[\alpha xt]}t \ . \]
We have \ \m{\phi(x)[\alpha xt]=x+\mu[\alpha xt]\alpha xt=x(1+\alpha\mu[\alpha 
xt]t)}. Hence
\[\epsilon' \ = \ \frac{\epsilon[\alpha xt]\alpha x}{\phi(x)[\alpha xt]} \ = \
\frac{\alpha\epsilon[\alpha xt]}{1+\alpha\mu[\alpha xt]t} \ \in R[n] \ . \]
We have \ \m{\epsilon'_0=\alpha\epsilon_0}, hence by proposition \ref{prop1}, 
if $\Theta$ is an isomorphism and $\alpha$ is invertible, \m{\epsilon_0} is 
invertible and so is \m{\epsilon'_0}, so \m{\Theta(\phi',\epsilon')} is an 
isomorphism.
\end{proof}

\sepprop

{\bf Remark: } suppose that $R$ is a $\C$-algebra. All the rings considered in 
\ref{PMR} are also $\C$-algebras, and we can restrict ourselves to morphisms of 
$\C$-algebras, and obtains the same results.

\end{sub}

\sepsub

\Ssect{Blowing-up of hypersurfaces}{BUH}

Let $X$ be a smooth, projective and irreducible variety. Let $\BS{\kx}$ be a
primitive multiple scheme of multiplicity \m{n\geq 2}, with underlying smooth 
variety $X$, and associated line bundle $L$ on $X$. Let \m{(U_i)_{i\in I}} be 
an affine open cover of $X$ such that we have trivializations
\xmat{\delta_i:U_i^{(n)}\ar[r]^-\simeq & U_i\times{\bf Z}_n , } (cf. 
\ref{PMS-1}). Let
\xmat{\delta_{ij}=\delta_j\delta_i^{-1}:U_{ij}\times{\bf Z}_n\ar[r]^-\simeq & 
U_{ij}\times{\bf Z}_n \ . }
Then \ \m{\delta_{ij}^*=\delta_i^{*-1}\delta_j^*} \ is an automorphism of \ 
\m{\ko_{U_i\times Z_n}=\ko_X(U_{ij})[t]/(t^n)}, such that for every \ 
\m{\phi\in\ko_X(U_{ij})[t]/(t^n)}, seen as a polynomial in $t$ with 
coefficients in \m{\ko_X(U_{ij})}, the term of degree zero of 
\m{\delta_{ij}^*(\phi)} is the same as the term of degree zero of $\phi$.

Let \m{H\subset X} be a hypersurface (i.e. a closed subscheme such that its 
ideal sheaf \m{I_H=I_{H,\BS{\kx}}} is invertible). We can choose the open cover 
\m{(U_i)_{i\in I}} such that, for every \m{i\in I}, \m{I_{H\cap U_i}} is 
generated by \m{x_i\in\ko_X(U_i)}. We have \ \m{{\bf B}_{H,X}=X}, and \m{{\bf 
B}_{H,\BS{\kx}}} is a primitive multiple scheme with underlying smooth variety 
$X$ 
and associated line bundle \m{L\ot\ko_X(H)}.

We now use \ref{PMR}, with \ \m{R=\ko_X(U_{ij})}. If \m{i,j\in I}, \m{i\not= 
j}, let
\[\theta_{ij} \ = \ \chi^K_{x_i}\circ\ov{\delta_{ij}^*}\circ\chi^K_{1/x_j} \ . 
\]
On \m{U_{ij}} we can write \ \m{x_i=\beta_{ij}.x_j}, with 
\m{\beta_{ij}\in\ko_X(U_{ij})} invertible. The family \m{(\beta_{ij})} is a 
cocycle defining the line bundle \m{\ko_X(H)}. By proposition \ref{prop2}, we 
have \ \m{\theta_{ij}=\ov{\tau}}, where \m{\tau} is an automorphism of 
\m{\ko_X(U_{ij})[t]/(t^n)}. We can write \ \m{\tau=\gamma_{ij}^*}, where 
\m{\gamma_{ij}} is an automorphism of \ \m{U_{ij}\times{\bf Z}_n}. The family 
\m{(\gamma_{ij}^*)} satisfies the cocycle relations \ 
\m{\gamma_{ij}^*\gamma_{jk}^*=\gamma_{ik}^*}, and defines a primitive multiple 
scheme \m{\BS{\kx}'} such that \ \m{\BS{\kx}'_{red}=X}. By formula 
$(\ref{equ1})$, \m{L(H)} 
is the line bundle on $X$ associated to \m{\BS{\kx}'}.

We have commutative diagrams
\xmat{\ko_X(U_{ij})[t]/(t^n)\ar[rr]^-{\delta_{ij}^*}\ar[d]^{\chi_{x_j}} & &
 \ko_X(U_{ij})[t]/(t^n)\ar[d]^{\chi_{x_i}} \\
\ko_X(U_{ij})[t]/(t^n)\ar[rr]^-{\gamma_{ij}^*} & & \ko_X(U_{ij})[t]/(t^n)
}
which define a morphism \ \m{\pi:\BS{\kx}'\to\BS{\kx}} \ which the identity on 
$X$, and 
such that its restriction to the open subset of $\BS{\kx}$ corresponding to 
\m{X\backslash H} is an isomorphism to the similar open subset of \m{\BS{\kx}'}.

%\sepprop
\newpage

\begin{subsub}\label{prop3}{\bf Proposition: } The morphism 
\m{\pi:\BS{\kx}'\to\BS{\kx}} is 
the blowing-up of $\BS{\kx}$ along $H$.
\end{subsub}
\begin{proof} We must prove that $\pi$ satisfies the following universal 
property: let $Y$ be a scheme and \ \m{f:Y\to\BS{\kx}} \ a morphism such that \
\m{\ki=f^{-1}(I_{H,\BS{\kx}}).\ko_Y} \ is an invertible sheaf of ideals on $Y$. 
Then 
there exists a unique morphism \ \m{\wt{f}:Y\to\BS{\kx}'} \ such that \ 
\m{\pi\circ\wt{f}=f}.

Let \m{y\in Y} be a closed point in the support of $\ki$ 
and suppose that \m{f(y)\in U_i}. Let \m{U\in Y} be an affine open neighborhood 
of $y$ such that \m{f(U)\subset U_i} and that \m{\ki_{|U}} is generated by 
\m{z\in\ki(U)}. It suffices to prove the claim for \m{f_{|U}}, since by
unicity, all the obtained morphisms \m{U\to\BS{\kx}'} can be glued to define 
\m{\wt{f}}. So we can assume that $f$ is a morphism \m{U\to U_i\times{\bf Z}_n}.

The morphism \ \m{\pi^*:\ko_X(U_i)[t]/(t^n)\to\ko_X(U_i)[t]/(t^n)} \ is defined 
by: for every \Nligne \m{a\in\ko_X(U_i)}, \m{\pi^*(a)=a}, and 
\m{\pi^*(t)=x_it}. We have \ \m{I_{H,U_i\times{\bf Z}_n}=(x_i,t)}. Now consider 
\Nligne \m{f^*:\ko_X(U_i)[t]/(t^n)\to\ko_Y(U)}. We can write \ 
\m{f^*(x_i)=\alpha z}, \m{f^*(t)=\beta z}, with \m{\alpha,\beta\in\ko_Y(U)}.

We will show that $\alpha$ is invertible.
We have \ \m{f^*(t^n)=0=\beta^nz^n}, hence \ \m{\beta^n=0}. If $\alpha$ is not 
invertible, he ideal \m{(\alpha)} can be included in a maximal ideal 
\m{\mm\not=\ko_Y(U)} of \m{\ko_Y(U)}. Since $\mm$ is prime and \m{\beta^n=0}, 
we have \m{\beta\in\mm}. By the definition of $\ki$, there exist 
\m{\lambda,\mu\in\ko_Y(U)} such that \ \m{\lambda f^*(x_i)+\mu f^*(t)=z}, i.e. 
\m{\lambda\alpha+\mu\beta=1\in\mm}. This is impossible, so $\alpha$ is 
invertible.

It follows that there one and only one possible definition for \ 
\m{\wt{f}^*:\ko_X(U_i)[t]/(t^n)\to\ko_Y(U)} : \m{\wt{f}^*(a)=f^*(a)} \ for 
every 
\m{a\in\ko_X(U_i)}, and \m{\wt{f}^*(t)=\beta\alpha^{-1}t}.
\end{proof}

\end{sub}

\sepsub

\Ssect{The case of double schemes}{DBS_H}

We suppose that \m{n=2} and keep the notations of the proof of theorem 
\ref{theo1}. Let \m{i,j\in I}. From \ref{DBC}, \m{\lambda_{ij}^*} and 
\m{\gamma_{ij}^*} are automorphisms of \m{\ko_X(U_{ij})[t]/(t^2)}, and there 
exists derivations \m{D_{ij}}, \m{E_{ij}} of \m{\ko_X(U_{ij})} such that, for 
every \m{u\in\ko_X(U_{ij})} we have
\[\lambda_{ij}^*(u) \ = \ u+D_{ij}(u)t \ , \qquad \gamma_{ij}^*(u) \ = \ 
u+E_{ij}(u)t \ . \]
According to \ref{DBC}, the family \m{(D_{ij})} (resp. \m{(E_{ij})}) represents 
\ \m{A\in H^1(X,T_X\ot L)} (resp. \ \m{\A\in H^1(X,T_X\ot L\ot\ko_X(H))}), and 
\m{\C.A} (resp. \m{\C.\A}) is the element of \m{{\bf DS}(X,L)} (resp. \m{{\bf 
DS}(X,L\ot\ko_X(H))}) corresponding to $\BS{\kx}$ (resp \m{\BS{\kx}'={\bf 
B}_{H,\BS{\kx}}}).

Let \ \m{\phi_H:H^1(X,T_X\ot L)\lra H^1(X,T_X\ot L\ot\ko_X(H))} \ be the 
canonical map.

\sepprop

\begin{subsub}\label{prop8}{\bf Proposition: } We have \ \m{\A=\phi_H(A)}.
\end{subsub}
\begin{proof} This follows easily from \ref{cech3} and proposition \ref{prop2}.
\end{proof}

\end{sub}

%\sepsec
\newpage

\section{Primitive double structures on \texorpdfstring{\m{\wt{\P}_2}}{5}}

\Ssect{The main result}{MA}

Let \ \m{P\in\P_2}. Let \ \m{\pi:\wt{\P}_2\to\P_2} \ be the blowing-up of 
\m{\P_2} along $P$, and \ \m{\wt{P}=\pi^{-1}(P)} \ the exceptional divisor.

According to \cite{dr10}, there exists only one non trivial primitive double 
scheme $\X$ such that \ \m{\X_{red}=\P_2}, and its associated line bundle is 
\m{\ko_{\P_2}(-3)}.

Let $\Y$ be a primitive double scheme such that \ \m{\Y_{red}=\P_2},  and with 
associated line bundle \m{\ko_{\P_2}(m)}. We will give the list of all 
primitive double schemes \m{\wt{\Y}} such that
\begin{enumerate}
\item[--] $\wt{\Y}_{red}=\wt{\P}_2$ ,
\item[--] there exists a morphism \ \m{\pi_{\wt{\Y}}:\wt{\Y}\to\Y} \ inducing 
\ \m{\pi_{P,\P_2}:\wt{\P}_2\to \P_2}, and such that its restriction \ 
\m{\wt{\Y}\backslash\wt{P}\to\Y\backslash\{P\}} \ is an isomorphism.
\end{enumerate}
In this case, if $\L$ is the line bundle on \m{\wt{\P}_2} associated to 
$\wt{\Y}$, there is an integer $p$ such that \ 
\m{\L=\pi_{P,\P_2}^*(\ko_{\P_2}(m))\ot\ko_{\wt{\P}_2}(p\wt{P})} . We will prove

\sepprop

\begin{subsub}\label{theo6}{\bf Theorem: } We have \ \m{p\geq 0}.

{\bf 1 -- } If \m{p=0}, the schemes $\wt{\Y}$ are the blowing-ups of $\Y$ along 
the good subschemes $Z$ such that \m{Z_{red}=\{P\}}. If \m{m\geq -1}, all 
these schemes are isomorphic. If \m{m\leq -2} they are all distinct.

{\bf 2 -- } If \m{p=1}, these schemes form a family 
\m{(\wt{\Y}_{1,\alpha})_{\alpha\in\C}}, \m{\wt{\Y}_{1,0}} being the blowing-up 
of the schemes of {\bf 1} along $\wt{P}$, and also the blowing-up of $\Y$ along 
$P$. If \m{\Y=\X}, these schemes are all distinct. If $\Y$ is trivial, 
\m{\wt{\Y}_{1,0}} is trivial, and all the \m{\wt{\Y}_{1,\alpha}}, 
\m{\alpha\not=0} are isomorphic and non trivial.

{\bf 3 -- } If \m{p\geq 2}, there is only one such scheme \m{\wt{\Y}_p}. If 
\m{p=2}, \m{\wt{\Y}_p} is the blowing-up of the schemes of {\bf 2} along 
$\wt{P}$. If \m{p>2}, \m{\wt{\Y}_p} is the blowing-up of \m{\wt{\Y}_{p-1}} 
along $\wt{P}$.
\end{subsub}

\sepprop

According to \ref{pt}, if \m{m\leq -2}, the set of schemes $\wt{\Y}$ such that 
\m{p=0} has a natural structure of affine space with associated vector space \ 
\m{T_{\P_2,P}\ot\ko_{\P_2}(m)_P}.

We will give the proof only for $\X$, the non trivial double scheme. The other 
cases are similar and simpler.

\end{sub}

\sepsub

\Ssect{Definition and properties of the non projective double plane 
\m{\X}}{x2_def}

Let $V$ be a complex vector space of dimension 3, \m{\P_2=\P(V)} \  and \ 
\m{\P_2^*=\P(V^*)}. Let $Q$ be the universal quotient bundle on \m{\P_2^*}. We 
have an exact sequence
\[0\lra\ko_{\P_2^*}(-1)\lra\ko_{\P_2^*}\ot V^*\lra Q\lra 0 \ . \]
%For every \m{x\in\P_2^*}, \m{x=\C.\phi} with \m{\phi\in V^*}, we have \ 
%\m{Q_x^*=\ker(\phi)}.

We fix a basis of $V$, so that we can define points of $V$ or \m{\P_2} using 3 
coordinates. For \m{0\leq i\leq 2}, let \ \m{U_i=\{(x_0,x_1,x_2)\in\P_2 \ ; \ 
x_i\not=0\}}. We will take the indices in \m{\Z/3\Z}, so that 
\m{x_3=x_0,x_4=x_1,\cdots}.

\sepprop

\begin{subsub}\label{x2_2} The double projective plane \m{\X} (cf. 
\cite{dr10},8) -- \rm We use the notations of \ref{blo} and 
\ref{PMS_def}.
For \m{0\leq i\leq 2}, \m{\Omega_{\P_2}(U_i)} is generated by 
\m{\dsp d\Big(\frac{x_{i+1}}{x_i}\Big)}, \m{\dsp 
d\Big(\frac{x_{i+2}}{x_i}\Big)}. The corresponding dual basis of 
\m{T_{\P_2}(U_i)} is \m{\dsp\Big(\frac{\partial}{\partial(x_{i+1}/x_i)},
\frac{\partial}{\partial(x_{i+2}/x_i)}\Big)}. Then $\X$ is constructed as in 
\ref{PMS-1}, \ref{DBC}, using the open cover \m{(U_i)_{0\leq i\leq 2}}. We take
\[\alpha_{i,i+1} \ = \ \Big(\frac{x_i}{x_{i+1}}\Big)^3 \ , \quad D_{i,i+1} \ = \
\frac{x_i}{x_{i+1}}\frac{\partial}{\partial(x_{i+2}/x_i)} \ . \]

So we have isomorphisms \ \m{\delta_i:U_i^{(2)}\lra U_i\times{\bf Z}_2} (where 
\m{U_i^{(2)}} is the open subscheme of \m{\X} corresponding to \m{U_i}). Let 
\[\delta_{ij}=\delta_j\delta_i^{-1}:U_{ij}\times{\bf Z}_2\hfl{\simeq}{}
U_{ij}\times{\bf Z}_2 \ , \]
and \ \m{\delta_{ij}^*=\delta_i^{*-1}\delta_j^*} \ is the corresponding 
automorphism of \m{\ko_{\P_2}(U_i)[t]/(t^2)}.

%For every \m{\alpha\in\ko_{\P_2}(U_i)}, we have
%\begin{eqnarray*}
%\delta_{i,i+1}^*(\alpha) & = & 
%\alpha+\frac{x_i}{x_{i+1}}\frac{\partial}{\partial(x_{i+2}/x_i)}(\alpha).t \ , 
%\\ \delta_{i,i+1}^*(t) & = & \Big(\frac{x_i}{x_{i+1}}\Big)^3t \ .
%\end{eqnarray*}
The line bundle on \m{\P_2} associated to \m{\X} is \m{\ko_{\P_2}(-3)}.
\end{subsub}

\end{sub}

\sepsub

\Ssect{The blowing-up of a point of \m{\P_2}}{blo}

Let \ \m{P=(0,0,1)\in\P_2} and \m{\U=\P_2\backslash\{P\}}. Let \ 
\m{\pi:\wt{\P}_2\to\P_2} \ be the blowing-up of \m{\P_2} along $P$. Let 
\m{\P\subset\P_2^*} be the line whose points are the lines of \m{\P_2} that 
contain $P$. We have
\[\wt\P_2 \ = \ \{(E,Q)\in\P\times\P_2\ ; \ Q\in E\} \ , \]
and \ \m{\pi=\pi_{P,\P_2}:\wt{\P}_2\to\P_2} \ is the second projection. Let \ 
\m{\Psi:\wt\P_2\to\P} \ be the first one. It is a projective bundle: 
\m{\wt\P_2=\P(Q^*_{|\P})}. Let \m{\ko_{\wt{\P}_2}(1)} be the corresponding line 
bundle, and \ \m{\wt{P}=\pi^{-1}(P)} \ the exceptional line. We have (cf. 
\cite{ha}, V, 2)
\[\ko_{\wt\P_2}(\wt{P}) \ \simeq \ 
\Psi^*(\ko_{\P}(-1))\ot\ko_{\wt\P_2}(1) \ , \quad
\pi^*(\ko_{\P_2}(1)) \ \simeq \ \ko_{\wt\P_2}(1) \ . \]

\sepprop

\begin{subsub}\label{op_co} An open cover -- \rm We have an isomorphism \ 
\m{D:\P_1\to\P} \ that associates to \m{(\alpha_0,\alpha_1)\in\P_1} the line of 
equation \m{\alpha_0x_0-\alpha_1x_1=0}. For \m{i=0,1}, let \
\m{V_i=\{(\alpha_0,\alpha_1)\in\P \ ; \ \alpha_i\not=0\}}. We have
\[\psi^{-1}(V_0) \ = \ \{(D(1,\lambda),(\lambda 
x_1,x_1,x_2))\in\wt\P_2 \ ; \ 
\lambda\in\C,(x_1,x_2)\in\P_1\} \ \simeq \ \C\times\P_1 \ , \]
\m{\psi^{-1}(V_0)\cap\U \ = \ U_1}, and a similar description of 
\m{\psi^{-1}(V_1)}.
Let \m{W_0=U_0}, \m{W_1=U_1},
\begin{eqnarray*}
W_2 & = & \{(E,(x_0,x_1,x_2))\in\psi^{-1}(V_0) \ ; \ x_2\not=0\} \\
& = & \{(D(1,\lambda),(\lambda\mu,\mu,1)) \ ; \lambda,\mu\in\C\} \ \simeq \  
\C^2 \ ,\\
W_3 & = & \{(E,(x_0,x_1,x_2))\in\psi^{-1}(V_1) \ ; \ x_2\not=0\} \\
& = & \{(D(\eta,1),(\nu,\eta\nu,1)) \ ; \eta,\nu\in\C\} \ \simeq \  
\C^2 \ .
\end{eqnarray*}
We have \ \m{W_0\cup W_1\cup W_2\cup W_3=\wt\P_2}, and
\[W_{02}=U_{012} \ , \ W_{03}=U_{02} \ , \ W_{12}=U_{12} \ , \ W_{13}=U_{012} 
\ , \ W_{23}\simeq \C\times\C^* \ , \]
\[W_{012} \ = W_{013} \ = \ W_{023} \ = \ W_{123} \ = \ U_{012} \ . \]
\end{subsub}

We will use the variables \ \m{\dsp\lambda=\frac{x_0}{x_1}}, \m{\dsp\mu=
\frac{x_1}{x_2}}. 

With the notations of \ref{pt}, we have \ \m{i_0=2}, \m{y_1=\mu}, 
\m{y_2=\nu=\lambda\mu}, \m{W_2=D^{\B_2}_+(y_1)} \ and \Nligne 
\m{W_3=D^{\B_2}_+(y_2)}.

\sepprop

\begin{subsub}\label{li_bu} Line bundles -- \rm With respect to the open cover 
\m{(W_i)_{0\leq i\leq 3}}, the line bundle \Nligne \m{\pi^*(\ko_{\P_2}(m))\ot
\ko_{\wt{\P}_2}(p\wt{P})} \ is represented by the cocycle 
\m{(\beta_{ij})_{0\leq i<j\leq 3}}, with
\[\beta_{01} \ = \ x_0^{-m}x_1^m \ = \ \lambda^{-m} \ , \quad \beta_{12} \ = 
\ x_1^{-p-m}x_2^{m+p} \ = \ \mu^{-p-m} \ , \]
\[\beta_{02} \ = \ x_0^{-m}x_1^{-p}x_2^{m+p} \ = \ \lambda^{-m}\mu^{-p-m} \ , 
\quad
\beta_{23} \ = \ x_0^{-p}x_1^{p} \ = \ \lambda^{-p} \ , \]
\[\beta_{03} \ = \ x_0^{-p-m}x_2^{m+p} \ = \lambda^{-p-m}\mu^{-p-m}
\ , \quad \beta_{13} \ = \ x_0^{-p}x_1^{-m}x_2^{m+p} \ = \ 
\lambda^{-p}\mu^{-p-m} \ 
. \]
\end{subsub}

\end{sub}

\sepsub

\Ssect{Morphisms to \m{\X}}{C-Y}

We use here the results of \ref{mor_db}.

\begin{subsub}\label{der} Rings and derivations -- \rm Let \m{\C(\lambda,\mu)} 
be the field of rational functions in the variables $\lambda$, $\mu$ (of 
\ref{op_co}). If \m{a_1,\ldots,a_n\in\C(\lambda,\mu)}, let 
\m{\BS{A}(a_1,\ldots,a_n)} be the sub-$\C$-algebra of \m{\C(\lambda,\mu)} 
generated by \m{a_1,\ldots,a_n}. We have
\[\ko_{\P_2}(U_2) \ = \ \BS{A}(\mu,\lambda\mu) \ , \quad
\ko_{\wt{\P}_2}(W_2) \ = \ \C[\lambda,\mu] \ , \quad
\ko_{\wt{\P}_2}(W_3) \ = \ \BS{A}(\frac{1}{\lambda},\lambda\mu) \ , \]
\[\ko_{\wt{\P}_2}(W_{12}) \ = \BS{A}(\lambda,\mu,\frac{1}{\mu})) \ , \quad
\ko_{\wt{\P}_2}(W_{23}) \ = \ \BS{A}(\lambda,\frac{1}{\lambda},\mu) \ , \
\ko_{\wt{\P}_2}(W_{03}) \ = \ \BS{A}(\frac{1}{\lambda},\lambda\mu,
\frac{1}{\lambda\mu}) \ , \]
\[\ko_{\P_2}(U_{012}) \ = \ \ko_{\wt{\P}_2}(W_{012}) \ = \
\BS{A}(\lambda,\frac{1}{\lambda},\mu,\frac{1}{\mu}) \ . \]

We will write the derivations $\kd$ of \m{\ko_{\wt{\P}_2}(W_i)} (or 
\m{\ko_{\wt{\P}_2}(W_{ij})}, \m{\ko_{\wt{\P}_2}(W_{012})}) as 
\[\kd \ = \ U.\frac{\partial}{\partial\lambda}+V.\frac{\partial}{\partial\mu} \ 
, \]
where \ \m{U,V\in\ko_{\wt{\P}_2}(W_{012})}. There will be conditions on $U$ and 
$V$, depending on the open subset where they are defined. For example suppose 
that we consider derivations of \m{\ko_{\wt{\P}_2}(W_{12})}. Then 
\m{\ko_{\P_2}(W_{12})=\BS{A}(\dsp\lambda,\mu,\frac{1}{\mu})}. So  $\kd$ is a 
derivation of \m{\ko_{\wt{\P}_2}(W_{12})} if and only if \ 
\m{\kd(\lambda),\kd(\mu)\in\ko_{\P_2}(W_{12})}, i.e. 
\m{U,V\in\ko_{\P_2}(W_{12})}. No verification is required for \m{W_{02}} and 
\m{W_{13}}, since \ 
\m{W_{02}=W_{13}=U_{012}}.
\end{subsub}

\sepprop

\begin{subsub} Construction of double structures -- \rm We want to define 
primitive double schemes $\Y$, such that \m{\Y_{red}=\wt{\P}_2}, and such that 
there exists a morphism \ \m{\phi:\Y\to\X} \ with \
\m{\phi_{red}=\pi:\wt{\P}_2\to\P_2}, and such that $\phi$ induces an 
isomorphism between the open subschemes of $\Y$, $\X$ over \ 
\m{\U=\P_2\backslash\{P\}}.

The associated line bundle will be of the form \ \m{L=\pi^*(\ko_{\P_2}(-3))\ot
\ko_{\wt{\P}_2}(p\wt{P})}, with associated cocycle \m{(\beta_{ij})} as in 
\ref{li_bu}, with \m{m=-3}.

We will use the open cover \m{(W_i)_{0\leq i\leq 3}} of \m{\wt{\P}_2} to build 
them as in \ref{PMS}. For this we need to define a derivation \m{E_{ij}} of 
\m{\ko_{\wt{\P}_2}(W_{ij}}), for \m{0\leq i<j\leq 3}, such that \ 
\m{(\beta_i^{-1}\ot E_{ij})_{0\leq i<j\leq 3}} represents the associated 
element of \m{H^1(\wt{\P}_2,T_{\wt{\P}_2}\ot L)}. In particular we have 
the relations \ \m{E_{ik}=E_{ij}+\beta_{ij}E_{jk}}. So we have, as in 
\ref{x2_def} for \m{\X}, for \m{0\leq i\leq 3}, isomorphisms \ 
\m{\rho_i:W_i^{(2)}\to W_i\times{\bf Z}_2}, with \ \m{\rho_0=\delta_0}, 
\m{\rho_1=\delta_1}.

We have \ \m{\dsp E_{01}=D_{01}=-\lambda^2\mu^2\frac{\partial}
{\partial\mu}}. We will take
\[E_{12} \ = \ A\frac{\partial}{\partial\lambda}+B\frac{\partial}
{\partial\mu} \ , \quad 
E_{13} \ = \ C\frac{\partial}{\partial\lambda}+D\frac{\partial}
{\partial\mu} \ . \]
Hence we have
\[E_{23} \ = \ \beta_{21}(E_{13}-E_{12}) \ = \ \mu^{p-3}\big(
(C-A)\frac{\partial}{\partial\lambda}+(D-B)\frac{\partial}{\partial\mu}\big)
\ , \]
\[E_{03} \ = \ E_{01}+\beta_{01}E_{13} \ = \ 
\lambda^3C\frac{\partial}{\partial\lambda}+(\lambda^3D-\lambda^2\mu^2)
\frac{\partial}{\partial\mu} \ . \]
The fact that \m{E_{ij}} is a derivation of \m{\ko_{\wt{\P}_2}(W_{ij})} will 
impose conditions on $A$, $B$, $C$ and $D$, necessary and sufficient to define 
$\Y$. The existence on $\phi$ will also impose other conditions.
\end{subsub}

\sepprop

\begin{subsub} Morphisms \m{\Y\to\X} -- \rm We have
\[\pi(W_0) \ = \ U_0 \ , \quad \pi(W_1) \ = \ U_1 \ , \quad \pi(W_2) \ 
\subset U_2 \ , \quad \pi(W_3) \ \subset \ U_2 \ , \]
and \ \m{\pi:W_i\to U_i} \ is the identity for \m{i=0,1}. Let
\[\Delta_0:W_0^{(2)}\lra U_0^{(2)} \ , \quad \Delta_1:W_1^{(2)}\lra U_1^{(2)}\]
be the identities. We will define 
\[\Delta_2:W_2^{(2)}\lra U_2^{(2)} \ , \quad \Delta_3:W_3^{(2)}\lra U_2^{(2)} \ 
, \]
(inducing $\pi$), such that \m{\Delta_1} and \m{\Delta_2} (resp. \m{\Delta_1} 
and \m{\Delta_3}) coincide on \m{W_{12}^{(2)}} (resp. \m{W_{13}^{(2)}}). It is 
easy to see that this describes completely a morphism \ 
\m{\phi:\Y\to\X}.

With the notations of \ref{mor_db3}, we have
\[\lambda_0=\delta_1:U_1^{(2)}\hfl{\simeq}{}U_1\times{\bf Z}_2 \ , \quad
\lambda_1=\delta_2:U_2^{(2)}\hfl{\simeq}{}U_2\times{\bf Z}_2 \ , \]
and \ \m{\lambda^*=\lambda_0^{*-1}\lambda_1^*=\delta_{12}^*} \ is an 
automorphism of \ \m{\ko_{\P_2}(U_{12})[t]/(t^2)}. We have
\[\lambda^*(a) \ = \ a+D_X(a)t \quad \text{for every \ } a\in\ko_{\P_2}(U_{12}) 
\ , \]
Where \ \m{\dsp D_X=\frac{x_1}{x_2}\frac{\partial}{\partial(x_0/x_1)}} \ (with 
respect to the variables \m{\dsp\frac{x_2}{x_1}}, \m{\dsp\frac{x_0}{x_1}}) and 
\ \m{\lambda^*(t)=\alpha t}, with \ \m{\dsp\alpha=\frac{x_1^3}{x_2^3}=\mu^3}. 
With respect to the variables \m{\lambda,\mu} we have \ 
\m{\dsp D_X=\mu\frac{\partial}{\partial\lambda}}.

We have similarly, with \ \m{\mu_0=\rho_1}, \m{\mu_1=\rho_2}, 
\m{\mu^*=\rho_{12}^*}, 
\[\mu^*(b) \ = \ b+D_Y(b)t \quad \text{for every \ } 
b\in\ko_{\wt{\P}_2}(W_{12}) \ , \]
where \m{D_Y} is the derivation \m{E_{12}} of the definition on \m{\Y}, and 
\ \m{\mu^*(t)=\beta t}, with \ 
\m{\dsp\beta=\beta_{12}=\mu^{-p+3}}. 

Let \ \m{A_2=\lambda_1\Delta_2\mu_1^{-1}=\delta_2\Delta_2\rho_2^{-1}:
W_2\times{\bf Z}_2\to U_2\times{\bf Z}_2}. We have
\[A_2^*(a) \ = \ \pi^*(a)+\gamma_2(a)t\]
for every \m{a\in\ko_X(U_2)}, where \ 
\m{\gamma_2:\ko_{\P_2}(U_2)\to\ko_{\wt{\P}_2}(W_2)} \ is a 
\m{\pi^*}-derivation, and \ \m{A_2^*(t)=\tau_2.t}, with 
\m{\tau_2\in\ko_{\wt{\P}_2}(W_2)}.

From $(\ref{equ30})$ we have \m{\dsp\tau_2=\mu^{p}}. It follows that
\m{p\geq 0}. From $(\ref{equ30})$ also we have
\[\gamma_2 \ = \ \mu^{p-3}\big((\mu-A)\frac{\partial}{\partial\lambda}-
B\frac{\partial}{\partial\mu}\big) \ . \]
We have \ \m{\ko_{\wt{\P}_2}(W_2)=\C[\lambda,\mu]} \ and \ 
\m{\ko_{\P_2}(U_2)=\BS{A}(\lambda\mu,\mu)}. So the fact that \m{\gamma_2} has 
values in \m{\ko_{\wt{\P}_2}(W_2)} is equivalent to \ 
\m{\gamma_2(\mu),\gamma_2(\lambda\mu)\in\C[\lambda,\mu]}. These conditions 
imply also that \m{E_{12}} is a derivation of \m{\ko_{\wt{\P}_2}(W_{12})}. 

We obtain similar results for \m{A_3}. We have
\[A_3^*(a) \ = \ \pi^*(a)+\gamma_3(a)t\]
for every \m{a\in\ko_X(U_2)}, where \ 
\m{\gamma_3:\ko_{\P_2}(U_2)\to\ko_{\wt{\P}_2}(W_3)} \ is a 
\m{\pi^*}-derivation, and \ \m{A_3^*(t)=\tau_3.t}, with 
\m{\tau_3\in\ko_{\wt{\P}_2}(W_3)}. From $(\ref{equ30})$ we have \ 
\m{\dsp\tau_3=\lambda^{p}\mu^{p}}. From $(\ref{equ30})$ we have
\[\gamma_3 \ = \ 
\lambda^{p}\mu^{p-3}\big((\mu-C)\frac{\partial}{\partial\lambda}-
D\frac{\partial}{\partial\mu}\big) \ . \]

\end{subsub}

\sepprop

\begin{subsub} Conditions on $A$, $B$, $C$, $D$ and consequences -- \rm We will 
study the conditions imposed on $A$, $B$, $C$ and $D$ by the facts that \ 
\m{\gamma_i(\ko_{\P_2}(U_2))\subset\ko_{\wt{\P}_2}(W_i)} \ for \m{i=2,3}, and 
that the \m{E_{ij}} are derivations of \m{\ko_{\wt{\P}_2}(W_{ij})}.

\end{subsub}

\sepprop

\begin{subsub}\label{lem3}{\bf Lemma: } We can suppose that 
\m{B=\mu^{-p+3}B_0}, where \m{B_0} is any element 
of \m{\C[\lambda,\mu]}, and that $A$ is of the form \ 
\m{A=\mu-\mu^{-p+2}R(\lambda)}, where $R$ is a polynomial.
\end{subsub}
\begin{proof}
We have \ \m{\mu\in\ko_{\P_2}(U_2)} \ and \ 
\m{\gamma_2(\mu)=-\mu^{p-3}B\in\ko_{\wt{\P}_2}(W_2)}, hence we can 
write \Nligne \m{B=\mu^{-p+3}B_1}, where \m{B_1\in\ko_{\wt{\P}_2}(W_2)}. But \ 
\m{\dsp B_1\frac{\partial}{\partial\mu}} \ is a derivation of 
\m{\ko_{\wt{\P}_2}(W_2)} and \ \m{\beta_{12}=\mu^{-p+3}}, hence by \ref{DBC} we 
can replace \m{E_{12}} with \ \m{\dsp 
E_{12}-\mu^{-p+3}(B_1-B_0)\frac{\partial}{\partial\mu}} , i.e. we can assume 
that \m{B=\mu^{-p+3}B_0}. 

We have \ \m{\gamma_2(\lambda\mu)=\mu^{p-2}(\mu-A)}, hence we can write \ 
\m{\mu^{p-2}(\mu-A)=P(\lambda,\mu)}, where $P$ is a polynomial, i.e. \ 
\m{A=\mu-\mu^{-p+2}P(\lambda,\mu)}. By \ref{DBC}, we can replace $A$ with 
\m{A+\mu^{-p+3}Q(\lambda,\mu)}, where $Q$ is a polynomial. The result for $A$ 
follows immediately.
\end{proof}

\sepprop

\begin{subsub}\label{lem4}{\bf Lemma: } Let \ \m{\dsp 
D=U\frac{\partial}{\partial\lambda}+V\frac{\partial}{\partial\mu}}, with
\m{U,V\in\C(\lambda,\mu)}. Then $D$ is a derivation of \m{\ko_{\wt{\P}_2}(W_3)} 
if and only there are \ \m{\epsilon,\tau\in\ko_{\wt{\P}_2}(W_3)} \ such that \ 
\m{U=\lambda^2\epsilon}, \m{V=\dsp\frac{1}{\lambda}\tau-\lambda\mu\epsilon}.
\end{subsub}
\begin{proof}
This follows easily from the fact that $D$ is a derivation of 
\m{\ko_{\wt{\P}_2}(W_3)} if and only if \ 
\m{D(1/\lambda),D(\lambda\mu)\in\ko_{\wt{\P}_2}(W_3)}.
\end{proof}

\sepprop

\begin{subsub}\label{lem5}{\bf Lemma: } We can suppose that
\[C \ = \ 
\mu+\lambda^{-p}\mu^{-p+2}\Big(-c_0\lambda+S\Big(\frac{1}{\lambda}\Big)
\Big) \ , \quad D \ = \ \lambda^{-p}\mu^{-p+3}c_0 \ , \]
where \m{c_0\in\C} and $S$ is a polynomial.
\end{subsub}
\begin{proof} We have \ 
\m{\gamma_3(\mu)=-\lambda^{p}\mu^{p-3}D\in\ko_{\wt{\P}_2}(W_3)} \ and \Nligne
\m{\gamma_3(\lambda\mu)=\lambda^{p}\mu^{p-3}(\mu(\mu-C))-\lambda D)
\in\ko_{\wt{\P}_2}(W_3)}. Hence
\[C \ = \ \mu-\lambda^{-p}\mu^{-p+3}\Big(-\frac{\lambda}{\mu}\gamma_3(\mu)+
\frac{1}{\mu}\gamma_3(\lambda\mu)\Big) \ ,
\quad D \ = \ -\lambda^{-p}\mu^{-p+3}\gamma_3(\mu) \ . \]
By \ref{DBC}, we can replace $C$ with \ 
\m{C+\lambda^{-p}\mu^{-p+3}.\lambda^2\epsilon} \ and $D$ with \ 
\m{\dsp D+\lambda^{-p}\mu^{-p+3}\Big(\frac{1}{\lambda}\tau-\lambda\mu\epsilon 
\Big)}, 
with \ \m{\epsilon,\tau\in\ko_{\wt{\P}_2}(W_3)}. We get
\[C \ = \ \mu+\lambda^{-p}\mu^{-p+3}\Big(\lambda^2\epsilon
+\frac{\lambda}{\mu}\gamma_3(\mu)-
\frac{1}{\mu}\gamma_3(\lambda\mu)\Big) \ ,
\quad D \ = \ 
\lambda^{-p}\mu^{-p+3}\Big(-\gamma_3(\mu)+\frac{1}{\lambda}\tau-
\lambda\mu\epsilon\Big) \ . \]
Since \ \m{\dsp\ko_{\wt{\P}_2}(W_3)=\BS{A}(\frac{1}{\lambda},\lambda\mu)}, we 
can take $\tau$, $\epsilon$ such that \ 
\m{-\dsp\gamma_3(\mu)+\frac{1}{\lambda}\tau-\lambda\mu\epsilon=c_0\in\C}. We 
have then \ \m{D=\lambda^{-p}\mu^{-p+3}c_0}. We can still replace $\epsilon$ 
with \ 
\m{\dsp\epsilon+\frac{1}{\lambda}\theta} \ and $\tau$ with \ 
\m{\tau+\lambda\mu\theta}, with \ \m{\theta\in\ko_{\wt{\P}_2}(W_3)}. We have
\[\lambda^2\epsilon
+\frac{\lambda}{\mu}\gamma_3(\mu)-\frac{1}{\mu}\gamma_3(\lambda\mu) \ = \
-\frac{\lambda}{\mu}c_0+\frac{1}{\mu}\big(\lambda\mu\theta+
\tau-\gamma_3(\lambda\mu)\big) \ . \]
We can then choose $\theta$ such that \ \m{\lambda\mu\theta+\tau-
\gamma_3(\lambda\mu)} \ is a polynomial in \m{\dsp\frac{1}{\lambda}}.
\end{proof}

\sepprop

Let \m{R_0} (resp. \m{S_0}) be the constant term of the polynomial 
\m{R(\lambda)} (resp. \m{\dsp S\Big(\frac{1}{\lambda}\big)}). Recall that we 
have \m{p\geq 0}.

\sepprop

\begin{subsub}\label{lem6}{\bf Lemma: } {\bf 1 -- } If \m{p\geq 2} then \ 
\m{R=S=c_0=0}.

{\bf 2 -- } If \m{p=1}, then \m{S=0} and \m{R=c_0}.

{\bf 3 -- } If \m{p=0}, then \m{R=R_0+c_0\lambda}, \m{S=S_0} and \m{R_0=-S_0}.
\end{subsub}
\begin{proof} We have \ \m{\dsp 
E_{23}(\lambda)=\mu^{p-3}(C-A)\in\ko_{\wt{\P}_2}(W_{23})=
\BS{A}(\lambda,\frac{1}{\lambda},\mu)}. The result follows immediately that, 
from lemmas \ref{lem3} and \ref{lem5}, we have
\[\mu^{p-3}(C-A) \ = \ \frac{1}{\mu}\Big(-c_0\lambda^{-p+1}+\lambda^{-p}S\Big(
\frac{1}{\lambda}\Big)+R(\lambda)\Big) \ . \]
\end{proof}

\sepprop

The fact that \m{E_{03}} is a derivation of \m{\ko_{\wt{\P}_2}(W_{03})} does 
not bring other conditions.
\end{sub}

\sepsub

\Ssect{proof of theorem \ref{theo6}}{pth6}

Suppose that \m{p=0}. Then we have, according to lemma \ref{lem6},
\[E_{12} \ = \ (\mu+\mu^2(-c_0\lambda-R_0))\frac{\partial}{\partial\lambda}+
\mu^3c_0\frac{\partial}{\partial\mu} \ , \]
and the same formula for \m{E_{13}}. With the notations and context of 
\ref{BL}, we have \m{i_0=2}, \ \m{y_1=\mu}, \m{y_2=\lambda\mu}, 
\m{J=\{0,1,i_1,i_2\}}, \m{(V_0,V_1,V_{i_1},V_{i_2})=(U_0,U_1,W_2,W_3)}. If \ 
\m{Z=Z(a_1,a_2)}, and \m{(E'_{ij})_{i,j\in J}} is the family of derivations 
that define \m{{\bf B}_{Z,\X}}, we have
\[E'_{1,i_1} \ = \ \mu\frac{\partial}{\partial\lambda}+\mu^3D(a_1,a_2) \ . \]
Using the preceding formulas, we have
\[\mu^3D(a_1,a_2) \ = \ \mu^3\big(a_1\frac{\partial}{\partial 
y_1}+a_2\frac{\partial}{\partial y_2}\big) \ = \ 
\mu^2(-a_1\lambda+a_2)\frac{\partial}{\partial\lambda}+\mu^3a_1\frac{\partial}
{\partial\mu} \ . \]
Hence we have \ \m{E_{12}=E'_{12}} \ if \ \m{a_1=c_0} \ \m{a_2=-R_0}. With 
simpler other checks, we see easily that the primitive double schemes of 
lemma \ref{lem6}, 3, are the \m{{\bf B}_{Z,\X}}.

The other cases (\m{p=1,2}) are similar. The scheme \m{\wt{\Y}_{1,\alpha}} of 
theorem \ref{theo6},2, corresponds to \m{c_0=\alpha} in lemma \ref{lem6}.

\end{sub}

\sepsub

\Ssect{The double schemes \m{\wt{\X}_{1,\alpha}}}{carpet}

We will study the schemes \m{\wt{\X}_{1,\alpha}}, \m{\alpha\in\C}, of theorem 
\ref{theo6}, 2 (when $\Y$ is the non trivial double scheme $\X$). The line 
bundle associated to \m{\wt{\X}_{1,\alpha}} is \ 
\m{\pi_{P,\P_2}^*(\ko_{\P_2}(-3))\ot\ko_{\wt{\P}_2}(\wt{P})=
\omega_{\wt{\P}_2}}. It follows that \m{\wt{\X}_{1,\alpha}} is a {\em 
K3-carpet} (cf. \cite{dr11}, \cite{ga-go-pu}).

\sepprop

\begin{subsub}\label{theo7}{\bf Theorem: } {\bf 1 -- } \m{\wt{\X}_{1,\alpha}} 
is 
quasi-projective if and only $\alpha\in\Q$ and \m{\alpha>0}.

{\bf 2 -- } If $\alpha\not\in\Q$, 
\m{\Pic(\wt{\X}_{1,\alpha})=\{\ko_{\wt{\X}_{1,\alpha}}\}}. 

{\bf 3 -- }If $\alpha\in\Q$, \m{\Pic(\wt{\X}_{1,\alpha})\simeq\Z}.
\end{subsub}

We will first describe \m{H^1(\wt{\P}_2,\Omega_{\wt{\P}_2})} with usual \v Cech 
cohomology, using the open cover \m{(W_i)_{0\leq i\leq 3}}. We have \ 
\m{H^1(\wt{\P}_2,\Omega_{\wt{\P}_2})\simeq\C^2}, and it is generated by 
\m{\nabla_0(\pi^*(\ko_{\P_2}(1)))} and\Nligne 
\m{\nabla_0(\ko_{\wt{\P}_2)}(-\wt{P})} (cf. \ref{can_class}).
With respect to \m{(W_i)_{0\leq i\leq 3}},
\begin{enumerate}
\item[--] $\nabla_0(\pi^*(\ko_{\P_2}(1)))$ is represented by the cocycle 
\m{(u_{ij})}, with
\[u_{01} \ = \ -\frac{d\lambda}{\lambda} \ , \quad u_{12} \ = \ 
-\frac{d\mu}{\mu} \ , \quad u_{23} \ = \ 0 \ . \]
\item[--] $\nabla_0(\ko_{\wt{\P}_2}(-\wt{P}))$ is represented by 
the cocycle \m{(v_{ij})}, with
\[v_{01} \ = \ 0 \ , \quad v_{12} \ = \ \frac{d\mu}{\mu} \ , \quad v_{23} \ = \ 
\frac{d\lambda}{\lambda} \ . \]
\end{enumerate}
For the canonical product \ \m{H^1(\wt{\P}_2,\Omega_{\wt{\P}_2})
\times H^1(\wt{\P}_2,\Omega_{\wt{\P}_2})\to 
H^2(\wt{\P}_2,\omega_{\wt{\P}_2})\simeq\C} , we have
\[\nabla_0(\pi^*(\ko_{\P_2}(1)))^2 \ = \ 1 \ , \quad 
\nabla_0(\ko_{\wt{\P}_2}(-\wt{P}))^2 \ = \ -1 \ , \quad
\nabla_0(\pi^*(\ko_{\P_2}(1))).\nabla_0(\ko_{\wt{\P}_2}(-\wt{P}))
\ = \ 0 \ . \]
From \ref{C-Y}, the scheme \m{\wt{\X}_{1,\alpha}} is defined, with respect to 
\m{(W_i)_{0\leq i\leq 3}}, by the family \m{(E_{ij})} of derivations, with
\[E_{01} \ = \ -\lambda^2\mu^2\frac{\partial}{\partial\mu} \ , \qquad
E_{12} \ = \ \mu(1-\alpha)\frac{\partial}{\partial\lambda} \ , \qquad
E_{23} \ = \ -\frac{\alpha}{\lambda}\frac{\partial}{\partial\mu} \ . \]
This family defines also \ \m{\sigma_\alpha\in 
H^1(\wt{\P}_2,\Omega_{\wt{\P}_2})} (which defines \m{\wt{\X}_{1,\alpha}} 
according to \ref{DBC}), in the sense of \ref{cech}.

\sepprop

\begin{subsub}\label{lem12}{\bf Lemma: } We have \ \m{\sigma_\alpha=
\nabla_0(\pi^*(\ko_{\P_2}(1)))+\alpha\nabla_0(\ko_{\wt{\P}_2}(-\wt{P}))} .
\end{subsub}
\begin{proof}
On the open subsets \m{W_i} we have the following trivializations \ 
\m{\alpha_i:\omega_{\wt{\P}_2|W_i}\to\ko_{W_i}} of \m{\wt{\P}_2} :
\[\alpha_0:d\Big(\frac{x_2}{x_0}\Big)\wedge d\Big(\frac{x_1}{x_0}\Big)\mapsto 1 
\ , \quad \alpha_1:d\Big(\frac{x_0}{x_1}\Big)\wedge 
d\Big(\frac{x_2}{x_1}\Big)\mapsto 1 \ , \quad \alpha_2:d\lambda\wedge 
d\mu\mapsto 1 \ , \quad \alpha_3:d\nu\wedge d\eta\mapsto 1 \ , \]
and, using the variables $\lambda$, $\mu$:
\[\alpha_0^{-1}:1\mapsto -\frac{1}{\lambda^3\mu^2}d\lambda\wedge d\mu \ , \ 
\alpha_1^{-1}:1\mapsto -\frac{1}{\mu^2}d\lambda\wedge d\mu \ , \
\alpha_2^{-1}:1\mapsto d\lambda\wedge d\mu \ , \
\alpha_3^{-1}:1\mapsto -\frac{1}{\lambda}d\lambda\wedge d\mu \ . \]
In usual \v Cech cohomology, \m{\sigma_\alpha} is defined by \ 
\m{E'_{ij}=\alpha_i^{-1}E_{ij}} :
\[E'_{01} \ = \ -\frac{1}{\lambda^3\mu^2}.d\lambda\wedge d\mu\ot(-\lambda^2\mu^2
\frac{\partial}{\partial\mu}) \ , \quad
E_{12} \ = \ -\frac{1}{\mu^2}.d\lambda\wedge 
d\mu\ot\mu(1-\alpha)\frac{\partial}{\partial\lambda} \ , \]
\[E_{23} \ = \ d\lambda\wedge 
d\mu\ot\big(-\frac{\alpha}{\lambda}\frac{\partial}{\partial\mu}\big) \ . \]
Now we apply the canonical isomorphism \ \m{\omega_{\wt{\P}_2}\ot
T_{\wt{\P}_2}\simeq\Omega_{\wt{\P}_2}} \ and we get
\[E_{01} \ = -\frac{1}{\lambda}.d\lambda \ , \quad
E_{12} \ = \ \frac{\alpha-1}{\mu}.d\mu \ , \quad
E_{23} \ = \ \frac{\alpha}{\lambda}.d\lambda \ . \]
The result follows from the preceding description of 
\m{\nabla_0(\pi^*(\ko_{\P_2}(1)))} and \m{\nabla_0(\ko_{\wt{\P}_2)}(-\wt{P})}.
\end{proof}

\sepprop

Now we prove theorem \ref{theo7}. A line bundle $\F$ on \m{\wt{\X}_{1,\alpha}} 
is ample if and only if \m{\F_{|\wt{\P}_2}} is ample (by \cite{ha2}, 
proposition 4.2). Hence \m{\wt{\X}_{1,\alpha}} is quasi-projective if and only 
if some ample line bundle on \m{\wt{\P}_2} can be extended to a line bundle to 
\m{\wt{\X}_{1,\alpha}}.

Let $F$ be a line bundle on \m{\wt{\P}_2}. According to \ref{ext_hig}, $F$ can 
be extended to a line bundle on \m{\wt{\X}_{1,\alpha}} if and only if \ 
\m{\sigma_\alpha.\nabla_0(F)=0}. Let \ 
\m{F_{mn}=\pi^*(\ko_{\P_2}(m))\ot\ko_{\wt{\P}_2}(-n\wt{P})}. Then we have \ 
\m{\sigma_\alpha.\nabla_0(F_{mn})=m-n\alpha}.

The ample line bundles on \m{\wt{\P}_2} are the \m{F_{mn}} such that \ 
\m{m,n>0}. It follows that \m{\wt{\X}_{1,\alpha}} is quasi-projective if and 
only if there exist positive integers $m$, $n$ such that \ \m{m-n\alpha=0}. 
This proves {\bf 1}.

According to \cite{dr11}, theorem 1.2.1, if \m{\F,\F'} are line bundles on 
\m{\wt{\X}_{1,\alpha}}, then \m{\F\simeq\F'} if and only if \ 
\m{\F_{|\wt{\P}_2}\simeq\F'_{|\wt{\P}_2}}. In other words, two extensions to  
\m{\wt{\X}_{1,\alpha}} of a line bundle on \m{\wt{\P}_2} are isomorphic. This 
implies {\bf 2} and {\bf 3}.

\end{sub}

\sepprop

\vskip 4cm

\end{document}